\documentclass{article}
\usepackage[utf8]{inputenc}
\usepackage{graphicx}
\usepackage[hyperref, dvipsnames]{xcolor}

\usepackage{amsmath, color, verbatim}
\usepackage[colorlinks=true]{hyperref}

\usepackage{amsfonts, mathtools, enumitem} 
\usepackage{hyperref}
\usepackage{todonotes}
\usepackage{enumitem}
\usepackage{titlesec}
 \usepackage{algorithm}
\usepackage{algorithmic}
 \usepackage{caption}
 \usepackage{placeins} 
\hypersetup{
  colorlinks=true,
}
\usepackage{amssymb}
\usepackage{amsthm}
\usepackage{bbm}
\usepackage[left=2cm,right=2cm,top=2cm,bottom=2cm]{geometry}
\usepackage{svg}

\usepackage[nameinlink,noabbrev,capitalize]{cleveref}

\DeclareMathOperator*{\esssup}{ess\,sup}

\DeclareMathOperator*{\dom}{dom}
\newcommand{\norm}[2]{\left\lVert #1\right\rVert_{#2}}

\newtheorem{thrm}{Theorem}[section]

\newtheorem{lmm}[thrm]{Lemma}
\newtheorem{prpstn}[thrm]{Proposition}

\newtheorem{ass}[thrm]{Assumption}

\newtheorem{crllr}[thrm]{Corollary}

\newtheorem{rmrk}[thrm]{Remark}
\newtheorem{xmpl}[thrm]{Example}

\newcommand{\grad}{\nabla}

\newcommand{\weaklyto}{\rightharpoonup}
\newcommand{\weakstar}{\stackrel{*}\rightharpoonup}
\newcommand{\cts}{\hookrightarrow}

\newcommand{\ctsCompact}{\xhookrightarrow{c}}
\newcommand{\E}{\mathbb{E}}
\newcommand{\pP}{\mathbb{P}}
\newcommand{\R}{\mathbb{R}}
\newcommand{\ER}{\overline{\mathbb{R}}}
\newcommand{\cR}{\mathcal{R}}
\newcommand{\cJ}{\mathcal{J}}

\newcommand{\CVaR}{\textup{CVaR}}
\newcommand{\D}{\textup{ d}}

\begin{document}
\hypersetup{
  urlcolor     = MidnightBlue, 
  linkcolor    = Bittersweet, 
  citecolor   = Cerulean
}
\title{Risk-averse optimal control of random elliptic variational inequalities
\footnotetext{We thank the anonymous referees for their careful reading and detailed comments, which helped to improve the paper.
AA and MH were partially supported by the DFG through the DFG SPP 1962 Priority Programme \emph{Non-smooth and Complementarity-based
Distributed Parameter Systems: Simulation and Hierarchical Optimization} within project 10. This research was conducted while the second author was affiliated with the Weierstrass Institute.
}}
\author{Amal Alphonse\thanks{Weierstrass Institute, Mohrenstrasse 39, 10117 Berlin, Germany (\href{mailto:alphonse@wias-berlin.de}{alphonse@wias-berlin.de})} 
\and Caroline Geiersbach\thanks{Fachbereich Mathematik, MIN Fakult\"at, Universit\"at Hamburg, Bundesstrrasse 55, 20146, Hamburg,
Germany; (\href{mailto:caroline.geiersbach@uni-hamburg.de}{caroline.geiersbach@uni-hamburg.de})} 
\and Michael Hinterm\"uller\thanks{Weierstrass Institute, Mohrenstrasse 39, 10117 Berlin, Germany (\href{mailto:hintermueller@wias-berlin.de}{hintermueller@wias-berlin.de})}
\and Thomas M.~Surowiec\thanks{Department of Scientific Computing and Numerical Analysis, Simula Research Laboratory, Kristian Augusts gate 23, 0164, Oslo, Norway (\href{mailto:thomasms@simula.no}{thomasms@simula.no})}}
\date{}

\newcommand{\skipeqs}[1]{\addtocounter{equation}{#1}}

\maketitle
\begin{abstract} We consider a risk-averse optimal control problem governed by an elliptic variational inequality (VI) subject to random inputs. By deriving KKT-type optimality conditions for a penalised and smoothed problem and studying convergence of the stationary points with respect to the penalisation parameter, we obtain two forms of stationarity conditions. The lack of regularity with respect to the uncertain parameters and complexities induced by the presence of the risk measure give rise to new challenges unique to the stochastic setting. We also propose a path-following stochastic approximation algorithm using variance reduction techniques and demonstrate the algorithm on a modified benchmark problem. \end{abstract}

\tableofcontents

\section*{Introduction}
In this work, we consider the following nonsmooth stochastic optimisation problem
\begin{equation}\label{eq:OCProblem}
\min_{u \in U_{ad}} \cR[\cJ(S(u))] +\varrho(u),
\end{equation}
where $U_{ad}$ is a set of controls, $S$ is the solution map of a random elliptic variational inequality (VI), 
$\cJ$ is an objective function, $\cR$ is a so-called risk measure that scalarises the random variable $\cJ(S(u))$, and $\varrho$ is the cost of the control $u$. 

Here, the state $S(u)=:y$ satisfies, on a pointwise almost sure (a.s.) level, the VI
\begin{equation}
\label{eq:VIProblemFormal}
\text{$y(\omega)  \leq \psi(\omega)$}: \quad \langle A(\omega)y(\omega) -f(\omega)-B(\omega)u, y(\omega) -v \rangle \leq 0 \quad \forall v : v \leq  \psi(\omega),
\end{equation} 
where $\omega \in \Omega$ stands for the uncertain parameter taken from a probability space $(\Omega, \mathcal{F}, \pP)$, $f(\omega)$ is a random source term,  $\psi(\omega)$ is a random obstacle, and $A(\omega)$ and $B(\omega)$ are random operators. The inequalities representing the obstacle constraint above should be understood in an almost everywhere sense, e.g., the first one is $y(\omega,x) \leq \psi(\omega, x)$ for a.e. $x \in D$ where $D$ is the spatial domain (to be specified below). The map $\cR$ is typically a convex functional chosen to generate solutions according to given risk preferences, e.g., optimal performance on average, weight of the tail of $\cJ(S(u))$, and so on.

Numerous free boundary problems in partial differential equations such as contact problems in mechanics and fluid flow through porous media  can be modelled as elliptic VIs, and in some cases, coefficients or inputs in the constitutive equations may be uncertain and are modelled as random. Such random VIs 
have been studied by \cite{Gwinner2000, Gwinner2006a, Forster2010, Kornhuber2014, Bierig2015}.  {
In addition, the mathematical programming literature offers an array of perspectives on variational inequalities subject to uncertainty \cite{Rockafellar2016}. These include the expected value (EV) approach \cite{Grkan1999}, the expected residual minimization (ERM) approach \cite{Chen2005}, and most recently \cite{Rockafellar2016,Rockafellar2018,Rockafellar2020}, which we will refer to as the Lagrangian (L) setting. These formulations differ from the random variational inequalities considered in this paper in several ways. The obvious differences are that the solution and probability spaces are typically finite dimensional/discrete, no differential operators are involved, and, especially in \cite{Rockafellar2016}, multi-stage formulations of the VI with a filtration of information are allowed. The single-stage EV, ERM, and L settings are closest to our random VI, however, the desired solution there needs to be deterministic, not a parametric response to the uncertainty as in our setting. As such, the overall goal is different: we seek a single, deterministic (here-and-now) decision $u$ that controls the distribution of values $X_u(\omega) := J(S(u))(\omega)$ according to the risk preference expressed by minimizing the risk measure $\mathcal{R}$ of $X_u$.}

It is important to already note here that \eqref{eq:OCProblem} typically contains two types of nonsmoothness: on the one hand from the solution operator $S$, on the other due to the choice of $\cR$, which in many interesting cases is a nonsmooth risk measure such as the {conditional} value-at-risk (usually written 
CVaR). The problem \eqref{eq:OCProblem}--\eqref{eq:VIProblemFormal} is formulated in the spirit of a ``here and now'' two-stage stochastic programming problem, where the decision $u$ is made before the realisation $\omega$ is made known. The study of such \textit{stochastic mathematical programs with equilibrium constraints} (SMPECs), has been limited to the finite-dimensional, risk-neutral (i.e., $\mathcal{R}=\E$) case; see \cite{Patriksson1999, Evgrafov2004a, shapiro2006stochastic, Shapiro2008}. 
For deterministic elliptic MPECs, there have been many developments in terms of theory and algorithms; see, e.g., \cite{Barbu1984, MR1739400, MR0423155, MR2515801, hintermuller2014several, MR2891941, MR2822818, surowiec2018numerical, MR3796767, Wachsmuth, hintermuller2008active, MR739836}.

The paper contains two contributions. First, using an adaptive smoothing approach, we derive stationarity conditions related to the well-known weak and C-stationarity conditions. To the best of our knowledge, this is the first attempt at such a derivation. {Second}, we provide a numerical study by applying a variance-reduced stochastic approximation method to solve an example of \eqref{eq:OCProblem}. 
Concerning the theoretical developments, our method for establishing the stationarity conditions is through a penalty approach similar to \cite{MR2822818, MR3056408}. 
As with the deterministic case, the penalty approach has the advantage that it is directly linked to the convergence analysis of solutions algorithms for the optimisation problem in a fully continuous, function space setting.  The theoretical results will also highlight a hidden difficulty unique to the stochastic setting that contrasts with the deterministic elliptic and parabolic cases. We also believe that this is an inherent difficulity in SMPECs in general, regardless of the dimension of the underlying decision space.

Our work is related to the problem setting in the recent paper~\cite{hertlein2022inexact} where the authors develop a bundle method for problems of the form \eqref{eq:OCProblem} with $\cR = \E$ the expectation.
The focus in our work is on obtaining stationarity conditions and a stochastic approximation algorithm for the general risk-averse case. Risk-averse optimisation is a subject in its own right; cf.~\cite{Shapiro2009} and \cite{Pflug2007} and the references therein. {Incorporating measures of risk in control problems appears to go back as far as \cite{Whittle1981}.} Modeling choices in engineering were explored in \cite{Rockafellar2015} and their application to PDE-constrained optimisation was popularised in \cite{Kouri2018, Kouri2019, Kouri2019a}. However, these papers typically require $S$ to be Fr\'echet differentiable, which does not hold in general for solution operators of variational inequalities. 
It is worth mentioning that typically
random VIs are studied in combination with some quantity of interest such as the expectation or variance as in \cite{Bierig2015, Kornhuber2014}. Another modeling choice could involve finding a deterministic solution $y$ satisfying \eqref{eq:VIProblemFormal}, which leads to an expected residual minimisation problem, or a $y$ satisfying the expected-value problem
\begin{equation*}
\text{$y  \leq \psi$}: \quad \langle \E[A(\cdot)y -f(\cdot)-B(\cdot)u], y -v \rangle \leq 0 \quad \forall v : v \leq  \psi.
\end{equation*}
These modeling choices are discussed in the survey \cite{Shanbhag2013}, but will not be pursued in the present work.

With respect to the organisation of the paper, after defining in Section \ref{sec:notation} some terms and notation,  Section \ref{sec:assumptions} is dedicated to introducing standing assumptions. An example of \eqref{eq:OCProblem} with specific choices for the various terms there is given in Section \ref{sec:example}. In Section \ref{sec:optimisation-problem}, the framework given above is formalised and it is shown that \eqref{eq:VIProblemFormal} has a unique solution (see  Lemma \ref{lem:boundOnSomega}). 
We will show in Proposition \ref{prop:propertiesOfS} that the control-to-state map $S$ maps into $L^q(\Omega;V)$ and hence the composition in \eqref{eq:OCProblem} is sensible. Moreover, we show in Proposition \ref{prop:existenceOC} that an optimal control to \eqref{eq:OCProblem} exists. In Section \ref{sec:penalised-control-problem}, we use a penalty approach on the obstacle problem and show that this penalisation is consistent with \eqref{eq:VIProblemFormal} (see Proposition \ref{lem:fundamentalLemma}). Optimality conditions for the control problem associated to the penalisation are given in Proposition \ref{prop:KKT-penalised} and Proposition \ref{prop:KKT-penalised-2} and form the starting point for the derivation of stationarity conditions. The main results culminate in Section \ref{sec:statConds}, namely Theorem \ref{thm:ss-weaker} and Proposition \ref{thm:ss}, where stationarity conditions of $\mathcal{E}$-almost weak and C-stationarity type are derived. This is done by taking the limit of the optimality conditions with respect to the penalisation parameter and is an especially delicate procedure due to the presence of the risk measure. A numerical example is shown in Section \ref{sec:numerics}. For the experiments, a novel path-following stochastic variance reduced gradient method is proposed in Algorithm \ref{alg:SVRG}. 
Although, as specified, we work in a particular setting of obstacle-type problems, we will explain in the concluding Section \ref{sec:conclusion} how greater generality could also be an option. 

\section{Preliminaries}
\subsection{Notation and background material}\label{sec:notation}

For exponents $t \geq 1$ {and given a Banach space $X$ and a probability space $(\Omega,\mathcal{F}, \pP)$, the Bochner space $L^t(\Omega; X):=L^t(\Omega,\mathcal{F}, \pP; X)$ is the set of all (equivalence classes of) strongly measurable functions $y\colon \Omega \rightarrow X$ having finite norm, where the norm is given by
\begin{equation*}
\lVert y \rVert_{L^t(\Omega;X)}:= \begin{cases}
                                     (\int_\Omega \lVert y(\omega) \rVert_X^t \D \pP(\omega))^{1/t}  &\text{ for } t \in [1, \infty),\\
                                    \esssup_{\omega \in \Omega} \lVert y(\omega) \rVert_X &\text{ for } t=\infty.
                                    \end{cases}
\end{equation*}}
We set $L^t(\Omega)=L^t(\Omega;\R)$ for the space of random variables with finite $t$-moments.  For a random variable $Z\colon \Omega \rightarrow \R$, the expectation is defined by $\E[Z]:=\int_{\Omega} Z(\omega)\D \pP(\omega).$
Recall that separability of {a reflexive Banach space implies separability of its dual space}. Moreover, if $X$ is a separable Banach space, then strong and weak measurability of the mapping $y\colon\Omega \rightarrow X$ coincide (cf.~\cite[Corollary 2, p.~73]{Hille1996}\footnote{As noted in \cite{Hille1996}, this result goes back to Pettis \cite{MR1501970} from 1938.}). Hence, we can call the mapping measurable without distinguishing between the associated concepts. Given Banach spaces $X,Y$, {denoting by $\mathcal{L}(X,Y)$  the set of bounded and linear maps from $X$ into $Y$}, an operator-valued function $A\colon \Omega \to \mathcal{L}(X,Y)$ is said to be \textit{uniformly measurable} in $\mathcal{F}$ if there exists a sequence of countably-valued operator random variables in $\mathcal{L}(X,Y)$ converging almost everywhere to $A$ in the uniform operator topology. 
A set-valued map {$T\colon \Omega \rightrightarrows X$} with closed images is called \textit{measurable} if the inverse image\footnote{For a set-valued map $T \colon\Omega \rightrightarrows X$ from $\Omega$ to a separable Banach space $X$,  the \textit{inverse image} on a set $E \subset X$ is
\begin{equation*}
    T^{-1}(E):= \{ \omega \in \Omega \colon T(\omega)\cap E \neq \emptyset \}.
\end{equation*}} of each open set is a measurable set, i.e., $T^{-1}(E) \in \mathcal{F}$ for every open set $E \subset X$.

Set $\ER := \mathbb{R} \cup \{\infty\}$. Recall that $F\colon X \to \ER$ is \textit{proper} if its \textit{effective domain}
\[\dom(F) := \{ x \in X : F(x) < \infty\}\]
satisfies $\dom(F) \neq \emptyset.$ 
The \textit{subdifferential} of a convex function $F\colon X \to \mathbb R$ at $z$ is the set $\partial F(z) \subset X^*$ defined as
\[\partial F(z) := \{ g \in X^* : F(z)-F(x) \leq \langle g, z-x\rangle_{X^*,X} \quad \forall x \in X\}.\]
We list some other notation and conventions that will be frequently used:
\begin{itemize}[leftmargin=*]\itemsep0em 
\item Whenever we write the duality pairing $\langle \cdot, \cdot \rangle$ without specifying the spaces, {we mean the pairing between the space $V$ (shortly to be introduced in \S \ref{sec:assumptions} below) and $V^*$}, i.e., $\langle \cdot, \cdot \rangle_{V^*,V}$.
\item For strong, weak, and weak* convergence, we use the symbols $\rightarrow$, $\rightharpoonup$ and $\weakstar$, respectively.
\item For a constant $t \in [1,\infty)$, $t'$ will denote its H\"older conjugate, i.e., $\tfrac{1}{t}+ \tfrac{1}{t'} = 1$.

\item We write $\cts$ to mean a continuous embedding and $\ctsCompact$ for a compact embedding.

\item $\mathcal{M}(\Omega;X)$ denotes the set of all {strongly $\mathcal{F}$-}measurable functions from $\Omega$ into $X$.

\item Statements that are true with probability one are said to hold almost surely (a.s.).

\item A generic positive constant that is independent of all other relevant quantities is denoted by $C$ and may have a different value at each appearance.

\end{itemize}
\subsection{Standing assumptions}\label{sec:assumptions}
Let us now describe our problem setup more precisely:
\begin{enumerate}[label=(\roman*)]\itemsep0em 
    \item $D \subset \mathbb{R}^d$ is a bounded Lipschitz domain for $d \leq 4$,  and take 
    \[H:=L^2(D) \quad\text{ and } V \in \{H^1(D), H^1_0(D)\}.\] 
    \item $U_{ad} \subset U$ is a non-empty, closed and convex set where the control space $U$ is a Hilbert space. 
    \item $(\Omega,\mathcal{F}, \pP)$ is a complete probability space, where $\Omega$ represents the sample space, $\mathcal{F} \subset 2^{\Omega}$ is the $\sigma$-algebra of events on the power set of $\Omega$, and $\pP\colon \Omega \rightarrow [0,1]$ is a probability measure.  
    \item $f \in L^r(\Omega;V^*)$ is a source term and $\psi \in L^s(\Omega;V)$ is an obstacle  for $r,s \in [2,\infty]$.
\item  
{$\cR\colon L^p(\Omega) \to \R$ and  $\varrho\colon U \to \R$} and  $p \in [1,\infty)$.
\end{enumerate}

To ease the presentation of the results, we make the following assumption on the almost everywhere boundedness of the operators in play in \eqref{eq:VIProblemFormal}. 
\begin{ass}\label{ass:newCombinedOnAandB}
The operators $A\colon \Omega \rightarrow \mathcal{L}(V,V^*)$ and $B\colon \Omega \rightarrow \mathcal{L}(U,V^*)$ are uniformly measurable and there exist positive constants $C_b, C_a, C_c$ such that for all $y, z \in V$ and $u\in U$ and a.s.~in $\omega \in \Omega$,
\begin{equation}\label{ass:boundedAndCoercive}
\begin{aligned}
     \langle A(\omega)y,y\rangle   &\geq C_a \lVert y\rVert_V^2,\\   
    \langle A(\omega)y, z \rangle   &\leq C_b\norm{y}{V}\norm{z}{V},\\
     \langle B(\omega)u, z \rangle  &\leq C_c\norm{u}{U}\norm{z}{V}.
\end{aligned} 
\end{equation}
\end{ass}
In applications, the operators $A$ and $B$ may be generated by random fields. There are numerous examples of random fields that are compactly supported; while this choice precludes a lognormal random field for $A$, in numerical simulations truncated Gaussian noise is often employed to generate samples. See, e.g.,  \cite{Lord2014, Ullmann2012, Gunzburger2014} for examples of compactly supported random fields, including approximations of lognormal fields as described.

The nature of the feasible set and composite objective function necessitates several assumptions on the objective functional. These ensure integrability, continuity and later, differentiability. 
\begin{ass}\label{ass:onJ}
Assume that $J\colon V \times \Omega \rightarrow \R$ is a Carath\'eodory {(that is, $J(v,\cdot)$ is measurable for {every} $v$ and $J(\cdot,\omega)$ is continuous for {$\pP$-a.e.}~$\omega.$)}  function and that there exists $C_1 \in L^p(\Omega)$ and $C_2 \geq 0$ such that
\begin{equation*}
    |J(v,\omega)| \leq C_1(\omega) + C_2\norm{v}{V}^{q\slash p},
\end{equation*}
where
 \begin{equation}
2 \leq q < \infty, \quad q \leq \min(r,s).  \label{eq:defnOfq}
\end{equation}
\end{ass}
For $y\colon \Omega \rightarrow V,$ we define the superposition operator $\mathcal{J}(y)\colon \Omega \rightarrow \R$ by $\cJ(y)(\omega) := J(y(\omega),\omega).$  The necessary and sufficient conditions to obtain continuity of $\cJ$ are directly related to famous results by Krasnosel’skii; see \cite{krasnoselskii} and \cite[Theorem 19.1]{vainberg1964variational}. Thanks to Assumption \ref{ass:onJ}, it follows by \cite[Theorem 4]{Goldberg} that 
\begin{equation}\label{eq:JisCts}
\cJ\colon L^q(\Omega;V) \to L^p(\Omega) \text{ is continuous}.
\end{equation}
\begin{rmrk}
{The conditions in \eqref{eq:defnOfq} force $q$ to be finite even if $f \in L^\infty(\Omega;V^*)$ and $\psi \in L^\infty(\Omega;V)$ (so that $r=s=\infty$)}. The case $q=\infty$ creates technical difficulties that we will address in a future work.
\end{rmrk}

In typical examples, $U \ctsCompact V^*$ and we would like the operator $B$ to mimic this compact embedding. {(Although $U$ is often taken to be $L^2(D)$ in the literature, other examples of $U$ one could consider include $\R^n$ or $L^2(\partial\Omega)$.)} For that purpose, we need the next assumption.
\begin{ass}\label{ass:BcompletelyCts}
If $u_n \weaklyto u$ in $U$ then $B(\omega)u_n \to B(\omega)u$ in $V^*$ a.s.
\end{ass}
Further assumptions will be introduced as and when required later in the paper. {We will take all of the above introduced assumptions as standing assumptions throughout the paper.}

\subsection{Example}\label{sec:example}
Take $V=H^1_0(D)$ and let $a\colon \Omega\times D \rightarrow \R$ be a given function such that $a_0 \leq a(\omega, x) \leq a_1$ a.s. and for a.e.~$x$, where $a_0>0$ and $a_1 > a_0$ are both constants. 
Define the operator
\[A(\omega) := -\grad \cdot (a(\omega)\grad u),\]
understood in the usual weak sense:
\[\langle A(\omega)y,z \rangle = \int_D a(\omega)\grad y \cdot \grad z \D x \qquad \text{for $y,z \in H^1_0(\Omega).$}\]
Set $U=L^2(D)$ with the box constraint set
\[U_{ad} := \{ u \in L^2(D) : u_a \leq u \leq u_b \text{ a.e.}\},\]
where $u_a, u_b \in L^2(D)$ are given functions. For $B$, we take it to be the canonical embedding $L^2(D) \ctsCompact H^{-1}(D)$, i.e., $B(\omega)u \equiv u$ as an element of $V^*=H^{-1}(\Omega)$. Take $f \in L^2(\Omega;H^{-1}(D))$, $\psi \in L^2(\Omega;H^1_0(D))$ and the exponent $q=2$. Let $p=1$  and define
\begin{equation}
\label{eq:tracking-type-Tikhonov}
J(y) := \frac 12\norm{y-y_d}{H}^2 \quad\text{and}\quad \varrho(u) := \frac{\nu}{2}\norm{u}{H}^2
\end{equation}
where $y_d \in L^2(D)$ is a given target state and $\nu > 0$ is the control cost. 
The risk measure is chosen to be the conditional value-at-risk, which for $\beta \in [0,1)$  is defined for a random variable $X\colon \Omega \rightarrow \R$ by
\begin{equation}
\label{eq:CVAR-def}
\mathcal{R}[X] = \CVaR_\beta[X] = \inf_{s \in \R} \left\lbrace s+ \frac{1}{1-\beta} \E[\max(X-s,0)] \right\rbrace.
\end{equation}
This risk measure is finite, convex, monotone, continuous and subdifferentiable {if $X \in L^1(\Omega)$} (see \cite[\S 6.2.4]{Shapiro2009}) and turns out to satisfy every assumption we will make in this paper. CVaR is easily interpretable: given a random variable $X$, CVaR$_{\beta}[X]$ gives the average of the tail of values $X$ beyond the upper $\beta$-quantile. The minimisers in \eqref{eq:CVAR-def} correspond to the $\beta$-quantile. $\CVaR_\beta$ approaches the essential supremum as $\beta \to 1.$

Further examples of risk measures can be found in \cite[\S 2.4]{Kouri2018} and references therein.

\section{Analysis of the optimisation problem}
\label{sec:optimisation-problem}
We begin by studying various properties of the solution map to the VI \eqref{eq:VIProblemFormal}, {which we restate here more precisely
\begin{equation*}
\text{$y(\omega) \in V, \; y(\omega)  \leq \psi(\omega)$}: \quad \langle A(\omega)y(\omega) -f(\omega)-B(\omega)u, y(\omega) -v \rangle \leq 0 \quad \forall v \in V: v \leq  \psi(\omega)\tag{\ref{eq:VIProblemFormal}}.
\end{equation*}
}
The solution mapping $u \mapsto y(\omega)$ in \eqref{eq:VIProblemFormal} is denoted by $S_\omega \colon U \to V$ and its associated superposition operator $S$ by
\begin{align}
    S(u)(\omega) := S_\omega(u).    \label{eq:defnOfSuperpositionS}
\end{align}
{We will later address the control problem \eqref{eq:OCProblem}.}
\subsection{Analysis of the VI}
We make heavy use, in particular, of the standing assumptions Assumption \ref{ass:newCombinedOnAandB}.

\begin{lmm}\label{lem:boundOnSomega}
For almost every $\omega \in \Omega$, there exists a unique solution to \eqref{eq:VIProblemFormal} satisfying the estimate
    \begin{align}
        \norm{S_\omega(u)}{V} &\leq  C\left(  \norm{f(\omega)}{V^*} + \norm{u}{U} + \norm{\psi(\omega)}{V}\right)\label{eq:solnMapPwBounded}
    \end{align}
where the constant $C>0$ depends only on $C_b, C_a$ and $C_c$. 
\end{lmm}
\begin{proof}
The conditions \eqref{ass:boundedAndCoercive} ensure the existence and uniqueness of the solution to \eqref{eq:VIProblemFormal} for {almost every} $\omega$ by the Lions--Stampacchia theorem; see \cite{lions1967variational}.

For the estimate we argue as follows. Setting $v=\psi(\omega)$ in \eqref{eq:VIProblemFormal} and splitting with Young's inequality {(i.e., the inequality $ab \leq \epsilon a^2 + b^2\slash (4\epsilon)$ for $\epsilon > 0$; we choose $\epsilon=C_a\slash 3$ to deal with the first two terms in the second line below)}, we obtain 
    \begin{align*}
        C_a\norm{y(\omega)}{V}^2 
        &\leq \langle A(\omega)y(\omega), \psi(\omega)\rangle +  \langle f(\omega)+B(\omega)u, y(\omega)-\psi(\omega) \rangle\\
        &\leq C_b\norm{y(\omega)}{V}\norm{\psi(\omega)}{V} + (\norm{f(\omega)}{V^*} + C_c\norm{u}{U})\norm{y(\omega)}{V}\\
        &\quad+ (\norm{f(\omega)}{V^*}+ C_c\norm{u}{U})\norm{\psi(\omega)}{V}\\
        &\leq \frac{C_a}{3}\norm{y(\omega)}{V}^2 + \frac{3C_b^2}{4C_a}\norm{\psi(\omega)}{V}^2 + \frac{3}{4C_a}(\norm{f(\omega)}{V^*} + C_c\norm{u}{U})^2 + \frac{C_a}{3}\norm{y(\omega)}{V}^2\\
        &\quad + \frac 12 (\norm{f(\omega)}{V^*} + C_c\norm{u}{U})^2  + \frac 12 \norm{\psi(\omega)}{V}^2.
    \end{align*}
This gives the uniform bound
\begin{align*}
\frac{C_a}{3}\norm{y(\omega)}{V}^2 &\leq \left(\frac{3C_b^2}{4C_a} + \frac 12\right)\norm{\psi(\omega)}{V}^2 + \left(\frac{3}{4C_a}+ \frac 12\right)(\norm{f(\omega)}{V^*} + C_c\norm{u}{U})^2.
\end{align*}
Thus we obtain \eqref{eq:solnMapPwBounded} if we take  
\[C \geq \max\left(\sqrt{\frac{9C_b^2}{4C_a^2} + \frac{3}{2C_a}}, \sqrt{\frac{9}{4C_a^2} + \frac{3}{2C_a}}, \sqrt{\frac{9C_b}{4C_a^2} + \frac{3C_b}{2C_a}}\right).\]
\end{proof}

\begin{rmrk}
In the proof above, we used as test function the obstacle $\psi$ since we assumed  (see Section \ref{sec:assumptions}) in particular that $\psi(\omega) \in V$ a.s.~and thus it is feasible. We could also have obtained a bound by testing with any map $v_0$ satisfying 
\begin{equation}\label{ass:feasibleElementR}
     v_0 \in \mathcal{M}(\Omega;V) : v_0(\omega) \leq \psi(\omega) \text{ a.s.}
 \end{equation}
 Assuming that such a map exists, this means that we could have asked for weaker regularity on $\psi$, e.g.~$\psi \in L^s(\Omega;L^2(D))$ with $\psi \geq 0$ on $\partial D$ would suffice (to the extent that the latter condition is defined). The proof of Proposition \ref{lem:boundsOnTtau} can be modified similarly too.
\end{rmrk}
\begin{rmrk}
If the constants $C_a, C_b, C_c$ in \eqref{ass:boundedAndCoercive} were functions defined on $\Omega$ instead,  we can prove an estimate similar to the one in Lemma \ref{lem:boundOnSomega}.  
For technical simplicity, we will not work in such generality in this paper. 
\end{rmrk}
For the next result and later ones, it is useful to note that, since $B(\omega)$ is bounded uniformly, Assumption \ref{ass:BcompletelyCts} on $B(\omega)$ being completely continuous also implies for $u_n \rightharpoonup u$ in $U$, by a simple Dominated Convergence Theorem (DCT) argument, that 
\begin{equation}
    Bu_n \to Bu \quad \text{in $L^t(\Omega;V^*)$ {for all $t \in [1, \infty)$}}.\label{eq:DCTConsequence}
\end{equation}
\begin{prpstn}\label{prop:propertiesOfS}
The solution map $S$ (defined in \eqref{eq:defnOfSuperpositionS}) satisfies the following.
\begin{enumerate}[label=(\roman*)]\itemsep0em 
    \item $S(u)\colon \Omega \to V$ is measurable for all $u \in U$.
    \item  $S(u) \in L^{\min(r,s)}(\Omega;V)$ for all $u \in U$.
    \item The estimate
    \[\norm{S_\omega(u) - S_\omega(v)}{V} 
    \leq C_a^{-1} \norm{B(\omega) (u-v)}{U}\] holds {a.s}.
    \item   
    If $u_n \weaklyto u$ in $U $, then 
    \[S_\omega(u_n) \to S_\omega(u)\text{ in $V$ a.s.}\]
    and
    \begin{align*}
    S(u_n)  &\to S(u) \text{ in $L^{q}(\Omega;V)$},\\
    S(u_n)-S(u)  &\to 0\text{ in $L^t(\Omega;V)$ {for all $t \in [1, \infty)$}}.
    \end{align*}
\end{enumerate}
\end{prpstn}
Note that the final claim is \textit{not} the same as $S(u_n) \to S(u) \text{ in $L^t(\Omega;V)$}$ because we do not know whether $S(u_n), S(u) \in L^t(\Omega;V)$ for {general} $t$.
\begin{proof}
\begin{enumerate}[label=(\roman*)]
    \item Let $u \in U$ be arbitrary but fixed and define the operator 
    \begin{equation*}
    {Q}\colon V \times \Omega \rightarrow V^*, (y,\omega) \mapsto Q(y,\omega):= A(\omega)y-f(\omega)-B(\omega)u.
    \end{equation*}
    The uniform measurability of $A$ and $B$ from Assumption \ref{ass:newCombinedOnAandB} implies (strong) measurability of $A(\cdot)y\colon \Omega \rightarrow V^*$ for every $y \in V$ and $B(\cdot)u\colon \Omega \rightarrow V^*$.  Thus ${Q}(y,\cdot)$ is measurable for every $y \in V$, and the almost sure continuity of ${Q}(\cdot,\omega)$ is clear. In particular, ${Q}$ is a Carath\'eodory  operator and is also superpositionally measurable \cite[Remark 3.4.2]{Gasinski2006}. 
    Note that the VI can be rewritten as: find $y(\omega) \leq \psi(\omega)$ such that $\langle {Q}(y(\omega),\omega), y(\omega)-v\rangle \leq 0$ for all $v \leq \psi(\omega)$.  
    Since the solution $y(\omega) = S(u)(\omega)$ exists by Lemma \ref{lem:boundOnSomega}, measurability follows from \cite[Theorem 2.3]{Gwinner2006a}. 
    \item This is an easy consequence of the estimate \eqref{eq:solnMapPwBounded}.    
    \item Let $y(\omega)=S_\omega(u)$ and $z(\omega)=S_\omega(v)$. Then for all $v \leq \psi(\omega)$ we have 
    \[ \langle A(\omega) (z(\omega)-y(\omega))-B(\omega)(v-u), z(\omega)-y(\omega)\rangle \leq 0.\]
    Using \eqref{ass:boundedAndCoercive},
    \[C_a \norm{y(\omega)-z(\omega)}{V}^2 \leq \norm{B(\omega)(u-v)}{V^*} \norm{z(\omega)-y(\omega)}{V}, \]
    and the claim follows.
    \item From the previous property, we derive
    \[\norm{S_\omega(u_n)-S_\omega(u)}{V} 
    \leq C_a^{-1}\norm{B(\omega)(u_n-u)}{V^*}.\]
    By Assumption \ref{ass:BcompletelyCts} we immediately obtain the pointwise a.s.~claim. Exponentating both sides to a power $t < \infty$, integrating and using \eqref{eq:DCTConsequence}, we obtain the Bochner convergence $S(u_n) -S(u) \to 0$ in $L^t(\Omega;V)$ {for all $t \in [1, \infty)$}, and by the sum rule, $S(u_n) \to S(u)$ in $L^{q}(\Omega;V)$ (recall that $q \leq \min(r,s)$ and $q \neq \infty$).\qedhere 
\end{enumerate}
\end{proof}
In the above result, we needed Assumption \ref{ass:BcompletelyCts} only for the final item (and later, it will be  needed only for the final item of Proposition \ref{lem:boundsOnTtau}).

\begin{rmrk}
It would also be possible to show measurability of the solution map if it is set-valued; see \cite{Gwinner2006a}. \end{rmrk}

\subsection{Existence of an optimal control}\label{sec:existenceOC}

We need certain basic structural properties on the functionals appearing in \eqref{eq:OCProblem} for the problem to be well posed. 
\begin{ass}\label{ass:onBothFunctionals}
Let {$\cR\colon L^p(\Omega) \to \R$} be lower semicontinuous and {$\varrho\colon U \to \R$} be  weakly lower semicontinuous. Assume also that
 \begin{equation}
 \text{either $U_{ad}$ is bounded or $\mathcal{R} \circ \cJ \circ S   + \varrho $ is coercive.}\label{ass:UadBdOrOpCoercive}    
 \end{equation}
 
\end{ass}
If a map from $L^p(\Omega)$ to $\mathbb{R}$ is finite, convex and monotone, it is continuous (and also subdifferentiable) on $L^p(\Omega)$ \cite[Proposition 6.5]{Shapiro2009}. 
\begin{prpstn}\label{prop:existenceOC}
Under 
Assumption \ref{ass:onBothFunctionals}, there exists an optimal control of \eqref{eq:OCProblem}.
\end{prpstn}
\begin{proof}
We can write the problem \eqref{eq:OCProblem} as
\[\min_{u \in U} F(u),\]
where $F\colon U \to \ER$ is defined by 
\[F := \cR \circ \cJ \circ S + \varrho + \delta_{U_{ad}},\]
with $\delta_{U_{ad}}\colon U_{ad} \to \{0, \infty\}$ denoting the indicator function, i.e., $\delta_{U_{ad}}(v)=0$ if $v\in U_{ad}$ and $\delta_{U_{ad}}(v)=\infty$ otherwise. We need to show that $F$ is proper, coercive and weakly lower semicontinuous.

The indicator function is proper since $U_{ad}$ is nonempty, hence 
$F$ is proper too. We have shown in Proposition \ref{prop:propertiesOfS} that $S\colon U \to L^q(\Omega;V)$ is completely continuous. This fact, combined with the continuity of $\cJ$ and lower semicontinuity of $\cR$ implies that $\cR \circ \cJ \circ S$ is weakly lower semicontinuous. Since $\varrho$ and $\delta_{U_{ad}}$ also possess this property (see, e.g., \cite[p.~10]{MR1727362} for the indicator function), we obtain weak lower semicontinuity of $F$.

Now, if $U_{ad}$ is bounded, $\delta_{U_{ad}}$ is coercive and we can use the fact that $\cR \circ \cJ \circ S + \varrho > -\infty$ to deduce the coercivity of $F$. Otherwise, if $\cR \circ \cJ \circ S + \varrho$ is coercive, we deduce the same property for $F$ by using the non-negativity of $\delta_{U_{ad}}.$

Applying now the direct method of the calculus of variations, 
we obtain the result.
\end{proof}

\begin{rmrk}
The conditions in Assumption \ref{ass:onBothFunctionals}, needed for the existence of controls, are  weaker than those needed in \cite[Proposition 3.12]{Kouri2018} (where  $\cR$ is taken to be finite, lower semicontinuous, convex
and monotone, and $\varrho$ to be proper, lower semicontinuous and convex) and weaker also than the ones in \cite[Proposition 3.1]{Kouri2019a} (where $\cR$ and $\varrho$ are finite, lower semicontinuous and convex and $\varrho$ is finite, convex and continuous). {A domain condition (of the same form as  \eqref{ass:domainCondition})} in the cited work is missing if $\varrho$ is taken to map into $\ER$, or otherwise, $\varrho$ needs to be finite, because  $\dom(\cR \circ \cJ \circ S + \varrho) \cap U_{ad}$ may be empty if this is not assumed. 
Note, however, that the authors above only have or assume weak-weak continuity for $S$ in their works.
\end{rmrk}

\section{A regularised control problem}
\label{sec:penalised-control-problem}

In this section, we approximate the VI \eqref{eq:VIProblemFormal} by a certain sequence of PDEs through a penalty approach. This leads to a regularised constraint in the overall control problem and is useful not only for numerical realisation and computation but also for the derivation of stationarity conditions for \eqref{eq:OCProblem}.
\subsection{A penalisation of the obstacle problem}
As alluded to, we penalise the constraint and essentially approximate \eqref{eq:VIProblemFormal} by a sequence of solutions of PDEs. These types of methods in the deterministic case have been comprehensively studied in the literature, see for example \cite[\S 3.5.2, p. 370]{Lions1969}, \cite[Chapter 1, \S 3.2]{GlowinskiLT} for some classical references.

Define the following differentiable approximation of the positive part  function $\max(0,\cdot)$:
\begin{equation}\label{eq:mrhoHK}
m_\tau(r) 
:= \begin{cases}
0 &: r \leq 0,\\
\frac{r^2}{2\tau} &: 0 < r < \tau,\\
r-\frac{\tau}{2} &: r \geq \tau.
\end{cases}
\end{equation}
We have $m_\tau \in C^1(\mathbb{R})$, $m_\tau' \in [0,1]$ and thanks to the regularity of its derivative when seen as function between the reals,  $m_\tau\colon H^1(D) \to L^2(D)$ is $C^1$ (see, e.g., \cite[Proposition 4]{MR2943603} for a direct proof). Before we proceed, note that other choices of penalisations or smooth approximations  $m_\tau$ are possible, see \cite{MR2822818, MR2891941, MR1739400}.

Consider for a fixed $\omega$ the penalised problem 
 \begin{equation}
 \label{eq:penalisedPDE}
 A(\omega)y +  \frac 1\tau m_\tau(y-\psi(\omega)) =f(\omega)+ B(\omega)u
 \end{equation}
 and denote the solution map $u \mapsto y$ as $T_{\tau,\omega} \colon U \to V$ and associated superposition map 
 $T_\tau(u)(\omega) := T_{\tau,\omega}(u).$  The solution map is well defined since \eqref{eq:penalisedPDE} has a unique solution \cite[Theorem 2.6]{Roubicek}. In an analogous way to the properties of $S$ obtained in Proposition \ref{prop:propertiesOfS}, we can show the following. 
\begin{prpstn}\label{lem:boundsOnTtau}
We have
\begin{enumerate}[label=(\roman*)]\itemsep0em 
    \item $T_\tau(u)\colon \Omega \to V$ is measurable for all $u \in U_{ad}$.
    \item The map {$T_\tau$ satisfies}
    \begin{align}
        &T_\tau \colon U \to L^{\min(r,s)}(\Omega;V)\label{eq:boundTtau}
    \end{align}
    and 
    \begin{equation}
    \label{eq:Ttau-bound}
    \norm{T_{\tau, \omega}(u)}{V} \leq  C\left(  \norm{f(\omega)}{V^*} + \norm{u}{U} + \norm{\psi(\omega)}{V}\right)
    \end{equation}
    holds {a.s}.
    \item The estimate 
    \begin{equation*}
        \norm{T_{\tau, \omega}(u)-T_{\tau, \omega}(v)}{V} \leq C_a^{-1}\norm{B(\omega)(u-v)}{V^*}
    \end{equation*}
    holds {a.s}.
    \item    
    If $u_n \weaklyto u$ in $U$, then 
    \[T_{\tau, \omega}(u_n) \to T_{\tau, \omega}(u) \text{ in $V$ a.s.},\]
    and 
    \begin{equation*}
    \begin{aligned}
        T_\tau(u_n) &\to T_\tau(u) &&\text{ in $L^{q}(\Omega;V)$},\\
        T_\tau(u_n) - T_\tau(u) &\to 0 &&\text{ in $L^t(\Omega;V)$ {for all $t \in [1, \infty)$}}.
    \end{aligned}
    \end{equation*}
\end{enumerate}
 \end{prpstn}
 \begin{proof}
 \begin{enumerate}[label=(\roman*)]
    \item This follows similarly to Proposition \ref{prop:propertiesOfS}.
     \item  Let $y=T_\tau(u)$ so that 
     \[ A(\omega)y(\omega) +  \frac 1\tau m_\tau(y(\omega)-\psi(\omega)) =f(\omega)+ B(\omega)u
\]
is satisfied almost surely. Testing the equation with  $y(\omega)-\psi(\omega)$ 
and using $m_\tau(r)r \geq 0$ leads to the same situation as in the proof of Lemma \ref{lem:boundOnSomega} and we obtain the same bound \eqref{eq:Ttau-bound}.
Then it follows as in Proposition \ref{prop:propertiesOfS} that $y$ is bounded uniformly in $L^{\min(r,s)}(\Omega;V)$ and the superposition map satisfies \eqref{eq:boundTtau}, just as for $S$.
\item If $y = T_{\tau, \omega}(u)$ and $z=T_{\tau, \omega}(v)$, we have
\begin{align*}
    \langle A(\omega)(z(\omega)-y(\omega)), z(\omega)-y(\omega) \rangle &+ \frac 1\tau\langle m_\tau(z(\omega)-\psi(\omega))-m_\tau(y(\omega)-\psi(\omega)), z(\omega)-y(\omega) \rangle\\
    &\quad= \langle B(\omega)(v - u), z(\omega)-y(\omega) \rangle,
\end{align*}
whence the estimate follows from the monotonicity of $m_\tau.$ 
\item From the previous property, we have
\begin{align*}
    \norm{T_{\tau, \omega}(u_n)-T_{\tau,\omega}(u)}{V} \leq C_a^{-1}\norm{B(\omega)(u_n-u)}{V^*}.
\end{align*}
From here, the argument is the same as in the proof of Proposition \ref{prop:propertiesOfS}. \qedhere
 \end{enumerate}

 \end{proof} 
The next result is fundamental as it shows that the solution of the penalised problem converges to the solution of the associated VI as the parameter is sent to zero. If we had defined $m_\tau \equiv m$ to be a penalty operator that is independent of $\tau$, we could apply classical approximation theory such as \cite[Theorem 5.2, \S 3.5.3, p.~371]{Lions1969}. The $m_\tau$-dependent case was addressed in \cite[Theorem 2.3]{MR2822818} but we need a slightly weaker assumption on the convergence of the source terms than was assumed there.

{In the following, all convergences involving $\tau$ correspond to the limit $\tau \downarrow 0$.}
\begin{prpstn}\label{lem:fundamentalLemma}
If $u_\tau \weaklyto u$ in $U$, then $T_{\tau, \omega}(u_\tau) \to S_\omega(u)$ in $V$ a.s. and $T_\tau(u_\tau) \to S(u)$ in $L^{q}(\Omega;V)$.
\end{prpstn}
\begin{proof}
The source term in the equation for $T_{\tau, \omega}(u_\tau)$ is $f(\omega) + B(\omega)u_\tau$ and by Assumption \ref{ass:BcompletelyCts}, this converges strongly to $f(\omega)+B(\omega)u$ in $V^*$ for almost every $\omega$. Therefore, we can apply \cite[Theorem 2.18]{AHROCQVI} (which yields the desired strong convergence  for a subsequence) and the subsequence principle to obtain $T_{\tau, \omega}(u_\tau) \to S_\omega(u)$ in $V$. By the bound \eqref{eq:Ttau-bound}, a simple DCT argument gives the Bochner convergence.
\end{proof}
Although \cite{Kouri2019a} provides a number of sufficient conditions for the continuity and differentiability properties of $T_{\tau}$ between $U$ and a corresponding Bochner space, not all cases are covered and \cite[Assumption 2.3]{Kouri2019a} appears to put a restriction on the integrability of the derivative of the nonlinearity that we may not have or need. In the following, we exploit the explicit structure of our nonlinearity and prove the necessary properties including continuous Fr\'echet differentiability directly.
 \begin{lmm}\label{lem:mTauIsC1}
 Let $\alpha \leq s$. The map $m_\tau \colon L^\alpha(\Omega;V) \to L^\beta(\Omega;H)$ is $C^1$ whenever $\beta < \alpha < \infty$. 
Furthermore, 
  \begin{equation*}
      \text{$m_\tau' \colon L^\alpha(\Omega;V) \to L^{\tilde\beta}(\Omega; \mathcal{L}(V, H))$ is continuous,}
  \end{equation*}
  where $\tilde\beta:=\alpha\beta\slash (\alpha-\beta)$.
 \end{lmm}
 \begin{proof}
 We want to apply the results in \cite{Goldberg}. 
 Define $G(\omega, y) := m_\tau(y-\psi(\omega))$ so that $G\colon \Omega \times V \to H$. 
 Set $\mathcal{G}(y)(\omega) := G(\omega, y(\omega))$ to be the associated Nemytskii operator. 
 
For $y \in V$,  we have 
 \begin{align*}
     \norm{G(\omega, y)}{H} 
     = \norm{m_\tau(y-\psi(\omega))}{L^2(D)} \leq  C\norm{y-\psi(\omega)}{L^2(D)},
 \end{align*}
 using the fact that $m_\tau$ is Lipschitz and $m_\tau(0)=0$. Hence 
 \begin{align*}
     \E \big[ \norm{\mathcal{G}(y)}{H}^\alpha\big] 
     \leq    C \E\big[ \norm{y-\psi(\cdot)}{L^2(D)}^\alpha \big],
 \end{align*}
which shows that $\mathcal{G}$  maps $L^\alpha(\Omega;V) \to L^\alpha(\Omega;H)$. Observe that we needed the assumption $\alpha \leq s$ for the right-hand side of the above to be finite. In particular, $\mathcal{G}\colon L^\alpha(\Omega;V) \to L^\beta(\Omega;H)$ whenever $\beta \leq \alpha.$ 

Now, as noted before, $m_\tau \colon V \to H$ is a $C^1$ map.  
Define an operator $K\colon \Omega \times V \to H$ by $K(\omega, y) := m_\tau'(y-\psi(\omega))$ and $\mathcal{K}(y)(\omega) := K(\omega, y(\omega))$.  We need $\mathcal{K}\colon L^\alpha(\Omega;V) \to L^{\tilde\beta}(\Omega; \mathcal{L}(V,H))$ to be continuous. 
If we had this, then applying \cite[Theorem 7]{Goldberg} we would obtain that $\mathcal G\colon L^\alpha(\Omega;V) \to L^\beta(\Omega;H)$ is a $C^1$ map for $\beta < \alpha$. We calculate
\begin{align*}
    \norm{m_\tau'(y-\psi(\omega))}{\mathcal{L}(V,H)} &= \sup_{\substack{g \in V\\ \norm{g}{V} = 1}}\norm{m_\tau'(y-\psi(\omega))g}{L^2(D)}\\
    &\leq  \sup_{\substack{g \in V\\ \norm{g}{V} = 1}}\norm{g}{L^2(D)}\\
    &\leq 1
\end{align*}
and hence, by \cite[Theorem 1]{Goldberg}, $\mathcal{K}$ maps all of $L^\alpha(\Omega;V)$ into $L^{\tilde\beta}(\Omega;\mathcal{L}(V,H))$. The continuity of $\mathcal{K}$ then follows from \cite[Theorem 4]{Goldberg}.
 \end{proof}

\begin{rmrk}
Using Lemma \ref{lem:mTauIsC1} and under the conditions on the exponents in the lemma, we can deduce that the map {$Q_\tau\colon U \times L^\alpha(\Omega;V) \to L^\beta(\Omega;V^*)$} is $C^1$, where
\[{Q_\tau(u,y) := Ay + \frac 1\tau m_\tau(y-\psi) -Bu - f}.\]
An application of the implicit function theorem would further restrict us to the setting where $\alpha \leq \beta$ in order to meet the condition of the isomorphism property for the partial derivative of {$Q_\tau$} and this is impossible. 
We will argue that $T_\tau$ is $C^1$ differently by directly using the equation it satisfies.
\end{rmrk}
 From now on, we need to explicitly use the fact that $V \cts L^4(D)$, which is true for $d \leq 4$ by Sobolev embeddings \cite[Theorem 4.12]{Adams2003}.

 \begin{prpstn}\label{lem:TtauFrechet}
The map $T_\tau\colon U \to L^q(\Omega;V)$ is continuously Fr\'echet differentiable and the derivative $T_\tau'(u)(h)$ belongs to $L^\infty(\Omega;V)$ and satisfies a.s.~the equation
   \begin{equation*}
 \label{eq:penalisedPDEDerivative}
 A(\omega)\delta(\omega) +  \frac 1\tau m_\tau'(y(\omega)-\psi(\omega))\delta(\omega) = B(\omega)h \qquad \text{where } y=T_\tau(u).
 \end{equation*}
For a given $h \in U$, the above equation admits a unique solution $\delta(\omega)\in V.$ Moreover, for all $\alpha \in [1,\infty]$ and $u \in U$,
 \[\frac{T_\tau(u+h)-T_\tau(u) - T_\tau'(u)(h)}{\norm{h}{U}} \to 0 \quad \text{in $L^\alpha(\Omega;V)$ as $h \to 0$ in $U$}.\]
  \end{prpstn}
It is worth emphasising that $T_\tau'(u) \in \mathcal{L}(U, L^\infty(\Omega;V))$ irrespective of $r$ and $s$ and that the quotient converges in $L^\infty(\Omega;V).$
 \begin{proof}
 Let $u, h \in U$. Then the existence and uniqueness of the solution $\delta$ follows from the monotonicity of $m_\tau'$. Define $y_h(\omega) := T_{\tau, \omega}(u+h)$, $y(\omega):=T_{\tau, \omega}(u)$ and the candidate derivative $\delta(\omega)$, which are defined in the following, on a pointwise a.s.~level:
 \begin{align*}
     A(\omega)y_h(\omega) + \frac 1\tau m_\tau(y_h(\omega)-\psi(\omega)) &= B(\omega)u + B(\omega)h + f(\omega),\\
     A(\omega)y(\omega) + \frac 1\tau m_\tau(y(\omega)-\psi(\omega)) &= B(\omega)u + f(\omega),\\
     A(\omega)\delta(\omega) + \frac 1\tau m_\tau'(y(\omega)-\psi(\omega))\delta(\omega) &= B(\omega)h.
 \end{align*}
Note that 
$\delta \in L^\infty(\Omega;V)$ due to Assumption \ref{ass:newCombinedOnAandB} and the non-negativity of $m_\tau'$. Let us for now omit the occurrences of $\omega$ for clarity. We have
 \begin{align*}
     A(y_h-y-\delta) + \frac 1\tau (m_\tau(y_h-\psi)-m_\tau(y-\psi)-m_\tau'(y-\psi)\delta) &= 0.
 \end{align*}
 Now, using the mean value theorem 
 on a pointwise a.e.~level in the domain $D$
 \begin{align*}
    & m_\tau(y_h-\psi)-m_\tau(y-\psi)-m_\tau'(y-\psi)\delta \\
     &\quad= \int_0^1 m_\tau'(y_h-\psi + \lambda(y-y_h))(y-y_h)\;\mathrm{d}\lambda  -m_\tau'(y-\psi)\delta\\
     &\quad= \int_0^1 m_\tau'(y_h-\psi + \lambda(y-y_h))(y-y_h)  -m_\tau'(y-\psi)\delta\;\mathrm{d}\lambda\\
     &\quad = \int_0^1 [m_\tau'(y_h-\psi + \lambda(y-y_h))-m_\tau'(y-\psi)](y-y_h)  + m_\tau'(y-\psi)(y_h-y-\delta)\;\mathrm{d}\lambda.
 \end{align*}
 Denote by $z:=y_h-y-\delta$. Plugging this in above, we find 
 \[Az + \frac{1}{\tau}\int_0^1 [m_\tau'(y_h-\psi + \lambda(y-y_h))-m_\tau'(y-\psi)](y-y_h)  + m_\tau'(y-\psi)z\;\mathrm{d}\lambda = 0.\]
 Testing with $z$, neglecting the final term due to non-negativity of $m_\tau'$, the above becomes
 \begin{align*}
 C\tau \norm{z}{V} &\leq 
 \left(\int_D \left|\int_0^1 [m_\tau'(y_h-\psi + \lambda(y-y_h))-m_\tau'(y-\psi)](y-y_h)\,\mathrm{d}\lambda \right|^2\mathrm{d}x\right)^{\frac 12}\\
 &\leq \frac 1\tau \left(\int_D \int_0^1 |y_h-y + \lambda(y-y_h)|^2|y-y_h|^2 \,\mathrm{d}\lambda\, \mathrm{d}x\right)^{\frac 12}\\
 &\leq \frac{\hat C}{\tau}\left(\int_D |y_h-y|^4 \,  \mathrm{d}x \right)^{\frac 12}\\
 &= \frac{\hat C}{\tau}\norm{y_h-y}{L^4(D)}^2\\
 &\leq \frac{\tilde C}{\tau} \norm{y_h-y}{V}^2
 \end{align*}
 using the embedding $V \cts L^4(D)$. Since $T_\tau$ is Lipschitz (see Proposition \ref{lem:boundsOnTtau}), we have
 \begin{align*}
 \frac{\tau^2 \norm{T_{\tau, \omega}(u+h)-T_{\tau, \omega}(u)-\delta(\omega)}{V}}{\norm{h}{U}}  &\leq C\norm{h}{U}.
 \end{align*}
 Taking the power $\alpha$, integrating over $\Omega$ and taking the $\alpha$-th root, we obtain the desired result for $\alpha < \infty$; taking the essential supremum recovers the remaining case. 
 \end{proof}
Let us now characterise the adjoint of the derivative of $T_{\tau,\omega}'$ (observe that $T_{\tau,\omega}'(u)^*\colon V^* \to U^*$ for $u \in U$). This will come in use later.
 \begin{lmm}\label{lem:characterisationOfAdjointTtauPrime}
For $u \in U$ and $g \in V^*$, we have
\begin{align*}
T_{\tau,\omega}'(u)^*g &= B(\omega)^*\eta(\omega)\quad\text{{a.s}},
\end{align*}
where $\eta(\omega) \in V$ satisfies {a.s}
\begin{align}
 A(\omega)^*\eta(\omega) + \frac 1\tau m_\tau'(y(\omega)-\psi(\omega))\eta(\omega) &= g, \quad\text{where $y(\omega) = T_{\tau,\omega}(u)$}.   \label{eq:foreta}
\end{align}
 If $g \in L^t(\Omega;V^*)$ {for some $t \in [1,\infty]$}, then $\eta \in L^t(\Omega;V)$.
\end{lmm}
\begin{proof}
Define pointwise a.s.~the quantity
 \[\lambda(\omega) := T_{\tau,\omega}'(u)^*g.\]
 Now, since $T_{\tau,\omega}'(u) = (A(\omega) + \frac 1\tau m_\tau'(y(\omega)-\psi(\omega)))^{-1}B(\omega),$ we find, using the commutation of adjoints and inverses, that 
 \[T_{\tau,\omega}'(u)^* = B(\omega)^*(A(\omega)^* + \frac 1\tau m_\tau'(y(\omega)-\psi(\omega)))^{-1}.\]
 Thus $\lambda(\omega)=B(\omega)^*(A(\omega)^* + \frac 1\tau m_\tau'(y(\omega)-\psi(\omega)))^{-1}g.$ If we then set $\eta(\omega):=(A(\omega)^* + \frac 1\tau m_\tau'(y(\omega)-\psi(\omega)))^{-1}g$, then
  \[\lambda(\omega)=B(\omega)^*\eta(\omega)\]
and $\eta(\omega)$ satisfies \eqref{eq:foreta}. 
 Testing the equation for $\eta$ with the solution itself and manipulating, we obtain
 \[C_a\norm{\eta(\omega)}{V} \leq \norm{g(\omega)}{V^*},\]
 from which we see that the integrability regularity is preserved.
\end{proof}

{
\begin{rmrk}
All of the results above would still hold if we had allowed for the possibility for $\cR$ and $\varrho$ to map into $\ER$, that is, if $\cR\colon L^p(\Omega) \to \ER$ and  $\varrho\colon U \to \ER$, as long as they are both proper and if we insert into Assumption \ref{ass:onBothFunctionals} the assumption that
\begin{equation}
    \dom( \cR \circ \cJ \circ S + \varrho) \cap U_{ad} \neq \emptyset.\label{ass:domainCondition}
\end{equation}
For reasons of clarity and simplicity, we have worked with finite risk measures and regularisers, and also because we need finiteness for what follows in Section \ref{sec:statRegOptControl}.
\end{rmrk}
}

 \subsection{Stationarity for the regularised control problem}\label{sec:statRegOptControl}

We consider the regularised problem
 \begin{equation}\label{eq:penalisedOCProblem}
 \min_{u \in U_{ad}} \cR[\cJ(T_\tau(u))] +\varrho(u)
 \end{equation}
 which, {if a condition akin to 
 \eqref{ass:UadBdOrOpCoercive} holds}, 
 by similar arguments to those made in  Section \ref{sec:existenceOC} has optimal controls $u_\tau^* \in U_{ad}$ with associated states 
 \[y_\tau^* := T_\tau(u_\tau^*) \in L^{\min(r,s)}(\Omega;V)\]
 (see Proposition \ref{lem:boundsOnTtau} for the integrability claim). The states satisfy
 \begin{equation}
      A(\omega)y_\tau^*(\omega) +  \frac 1\tau m_\tau(y_\tau^*(\omega)-\psi(\omega)) =f(\omega)+ B(\omega)u_\tau^*\quad \text{a.s.}\label{eq:ytau}
 \end{equation}
 We will show in Proposition \ref{prop:convergenceOfOpts} that the {optimal pairs} $(y_\tau^*, u_\tau^*)$ converge to {an optimal pair} of the original control problem \eqref{eq:OCProblem}.

We make the following standing assumption (in addition to the conditions assumed in Section \ref{sec:assumptions}). 
The assumption introduces requirements on $J_y$ that are needed to ensure the Fr\'echet differentiablity of $\cJ$. We have to slightly restrict the range of exponents that were available in Assumption \ref{ass:onJ} in order to avoid trivial situations (see the comment after Lemma \ref{lem:goldberg}). Note that if we had allowed for $q=\infty$, we would also need $p=\infty$ in order to apply the theory of \cite{Goldberg} to get $\cJ$ to be $C^1$, however $p=\infty$ is problematic, see Remark \ref{rem:infty}.

\begin{ass}\label{ass:onvarrho}
Assume the following.
\begin{enumerate}[label=(\roman*)]\itemsep0em 
\item\label{item:JandJyCC}  The function $J(\cdot,\omega)\colon V \to \mathbb{R}$ is a.s.~continuously Fr\'echet differentiable with $J_y(\cdot, \omega)\colon V \to V^*$ Carath\'eodory. 

\item\label{ass:finiteCase} Take 
\[p < q < \infty\]
and defining
\[\tilde p:= \frac{pq}{q-p},\] 
 there exists $\tilde{C}_1 \in L^{\tilde{p}}(\Omega)$ and $\tilde{C}_2 \geq 0$ such that
\begin{equation}
\norm{J_y(v,\omega)}{V^*} \leq \tilde{C}_1(\omega) + \tilde{C}_2\norm{v}{V}^{q/\tilde{p}}. 
\label{ass:JyBoundedFinite}
\end{equation}

\item\label{item:coerciveOrBddForTtauAndS} Suppose 
 \begin{equation*}
  \text{either $U_{ad}$ is bounded or both $\mathcal{R} \circ \cJ \circ S  + \varrho $ and $ \mathcal{R} \circ \cJ \circ T_\tau + \varrho $  are coercive.}
 \end{equation*}
\item\label{item:newOnR} The map  $\mathcal{R}\colon L^p(\Omega) \to \mathbb R$ is convex and lower semicontinuous. 
\item\label{item:varrhoConvexGateaux} The map $\varrho\colon U \to \mathbb{R}$ is convex and G\^{a}teaux directionally differentiable with $\varrho'(u)(v)$ denoting the G\^{a}teaux derivative at $u \in U$ in the direction $v \in U$.
\item\label{item:wlscOfVarrhoPrime} The map $\varrho'(\cdot)(\cdot-v)\colon U \to \mathbb{R}$ is weakly lower semicontinuous for every $v \in U_{ad}$.
\end{enumerate}
\end{ass}

Under items \ref{item:JandJyCC} and \ref{ass:finiteCase}, we obtain through Lemma \ref{lem:goldberg} that $\cJ\colon L^q(\Omega;V) \to L^p(\Omega)$ is continuously Fr\'echet differentiable with
\[\mathcal{J}'(y)(h)
= J_y(y(\cdot), \cdot)h(\cdot).\]
Furthermore, 
\begin{equation}
\text{$\mathcal{J}'\colon L^q(\Omega;V) \to L^{\tilde p}(\Omega;V^*)$ is continuous}\label{eq:JyCtsGoodSpace}.
\end{equation}
A few additional words on these assumptions are in order.
\begin{rmrk}\label{rem:secondSetOfAss}
\begin{enumerate}[label=(\roman*)]\itemsep0em 
\item The coercivity part of Assumption \ref{ass:onvarrho} \ref{item:coerciveOrBddForTtauAndS} is automatic in the typical case where $\cR$ is monotonic, $J$ is a tracking-type functional (thus $\cR\circ \cJ \circ T_\tau \geq 0$) and $\varrho(u) = (\nu\slash 2)\norm{u}{U}^2$ where $\nu>0$ is a constant.
\item By Assumption \ref{ass:onvarrho} \ref{item:newOnR}, it follows that $\cR$ is continuous (since it is finite, lower semicontinuous and convex), see \cite[Proposition 2.111]{Bonnans2013}. Then by \cite[Proposition 2.126 (v)]{Bonnans2013}, $\cR$ is Hadamard differentiable with
\[\cR'[z](h) = \sup_{\nu \in \partial\cR(z)} \E[\nu h ].\]
\item It is easy to see that
\begin{equation}
\tilde p > p\quad\text{and}\quad \tilde p \geq q'.\label{eq:tildePGreaterThanP}    
\end{equation}
Some further relations between the various exponents can be found in Section \ref{sec:pqStuff}.
\item Assumption \ref{ass:onvarrho} \ref{item:varrhoConvexGateaux} implies weak lower semicontinuity of $\varrho$. Indeed, take a sequence $w_n \weaklyto w$ in $U$. Since $\varrho$ is G\^ateaux differentiable and convex, we have
\[\varrho(w_n) \geq \varrho(w) + \varrho'(w)(w_n-w),\]
and we get the claim by using the fact that the derivative of $\varrho$ is linear in the direction.
\end{enumerate}
\end{rmrk}

\begin{xmpl}
If $J$ is of tracking type as in \eqref{eq:tracking-type-Tikhonov}, we have that $J_y(y, \omega) = (y-y_d, \cdot)_H$
so that
$\norm{J_y(y,\omega)}{V^*} \leq C(\norm{y_d}{H} + \norm{y}{V})$ and the growth condition assumption  \eqref{ass:JyBoundedFinite} is satisfied. 
\end{xmpl}

Now, we have that $\mathcal{J} \circ T_\tau\colon U \rightarrow L^p(\Omega)$ is $C^1$ since $T_\tau\colon U \to L^{q}(\Omega;V)$ and $\mathcal{J}\colon L^q(\Omega;V) \to L^p(\Omega)$ are $C^1$ (see Proposition \ref{lem:TtauFrechet} for the former). Since $\cR$ is Hadamard differentiable, the chain rule \cite[Proposition 2.47]{Bonnans2013} yields that 
\[(\cR \circ \cJ \circ T_\tau)'(u_\tau^*)(h) = \cR'(\cJ(y_\tau^*))\cJ'(y_\tau^*)T_\tau'(u_\tau^*)h \quad \forall h \in U.\]
Using the expression for the derivative $\cR'$ above and the fact that  $\cJ'(y_\tau^*)T_\tau'(u_\tau^*)\colon U \to L^p(\Omega)$ is continuous, arguing in the usual way for B-stationarity (see \cite{hintermuller2014several} for the corresponding notion),
we obtain
\begin{equation*}
    \sup_{\pi \in \partial \cR[\cJ(y_\tau^*) ]} \E[\cJ'(y_\tau^*)(T'_\tau(u_\tau^*)h)\pi] +\varrho'(u_\tau^*)(h) \geq 0 \quad \forall h \in T_{U_{ad}}(u_\tau^*),
\end{equation*}
where  $T_{U_{ad}}(u_\tau^*)$ is the tangent cone of $U_{ad}$ at $u_\tau^*.$ This condition is not so convenient due to the supremum present. However, making use of subdifferential calculus, we have the following characterisation.
\begin{lmm}
\label{lem:optimality-penalized}
{Let Assumption \ref{ass:onvarrho} hold.} There exists $\pi_\tau^* \in \partial \cR[\cJ(y_\tau^*)]$ such that
\begin{equation}
\label{eq:OC-penalized-as-VI} \E[\cJ'(y_\tau^*)(T'_\tau(u_\tau^*)(u_\tau^*-v))\pi_\tau^*] + \varrho'(u_\tau^*)(u_\tau^*-v) \leq 0 \quad \forall v \in U_{ad}.
\end{equation}
\end{lmm}
\begin{proof}
We begin by checking some properties of $\cR\circ \cJ \circ T_\tau$ that we need to apply the sum rule for subdifferentials.

Note first that 
$\cR$ is locally Lipschitz  on $L^p(\Omega)$ (which means that it is locally Lipschitz near every point of $L^p(\Omega)$). Also, $\mathcal{J} \circ T_\tau\colon U \rightarrow L^p(\Omega)$  is strictly differentiable  at $u_\tau^*$ (in the sense of Clarke) since it is $C^1$ (see \cite[Corollary, p.~32]{Clarke}), 
and by \cite[Proposition 2.2.1]{Clarke}, it is Lipschitz near $u_\tau^*$. It follows that the composition $\cR\circ \cJ \circ T_\tau$ is locally Lipschitz near $u_\tau^*$.

Since $\varrho$ is convex and weakly lower semicontinuous (see Remark \ref{rem:secondSetOfAss}), 
$\varrho$ is also locally Lipschitz near every point of $U$. Thus the sum $\cR\circ \cJ \circ T_\tau + \varrho$ is also locally Lipschitz near $u_\tau^*$. Hence by the corollary on page 52 of \cite{Clarke}, we have
\[0 \in \partial(\cR \circ \cJ \circ T_\tau + \varrho)(u_\tau^*) + N_{U_{ad}}(u_\tau^*)\]
where $N_{U_{ad}}(u)$ stands for the normal cone of $U_{ad}$ at $u$. By \cite[Proposition 2.3.3]{Clarke}, we get
\[ \partial(\cR \circ \cJ \circ T_\tau + \varrho)(u_\tau^*)  \subset  \partial(\cR \circ \cJ \circ T_\tau)(u_\tau^*) + \partial \varrho(u_\tau^*) .\]
With $\cR$ being locally Lipschitz on $L^p(\Omega)$, we can apply the subdifferential chain rule \cite[Theorem 2.3.10 and Remark 2.3.11]{Clarke} to obtain
\[\partial (\cR\circ \cJ \circ T_\tau)(u_\tau^*) = [(\cJ \circ T_\tau)'(u_\tau^*)]^*\partial \cR(\cJ \circ T_\tau(u_\tau^*)).\]
Equality here holds since $\cR$ is convex and thus regular \cite[Proposition 2.3.6 (b)]{Clarke}. Consolidating all of the above, we have
\[0 \in [(\cJ \circ T_\tau)'(u_\tau^*)]^*\partial \cR(\cJ \circ T_\tau(u_\tau^*)) + \partial \varrho(u_\tau^*) + N_{U_{ad}}(u_\tau^*).\]
Now we argue similarly to the corrigendum to \cite{Kouri2018}.  It follows that there exists a $\nu \in \partial \cR(\cJ(y_\tau^*))$ and $\eta \in \partial\varrho(u_\tau^*)$ satisfying
\[-[(\cJ \circ T_\tau)'(u_\tau^*)]^*\nu  - \eta  \in N_{U_{ad}}(u_\tau^*).\]
By \cite[Proposition 2.1.2 (b) and Proposition 2.2.7]{Clarke}, $\varrho'(u_\tau^*)$ is the support function of $\partial\varrho (u_\tau^*)$, i.e.,
\[\varrho'(u_\tau^*)(h) = \sup_{g \in \partial\varrho(u_\tau^*)} \langle g, h \rangle_{L^{p'}(\Omega),L^p(\Omega)} \quad \forall h \in U\]
and hence, by linearity of the derivative, $\eta=\varrho'(u_\tau^*)$. By definition of the normal cone and using the convexity of $U_{ad}$, we then have that $\nu \in \partial\cR(\cJ(y_\tau^*))$ satisfies
\[\langle [\cJ'(y_\tau^*)T_\tau'(u_\tau^*)]^*\nu + \varrho'(u_\tau^*), u_\tau^*- z \rangle_{U^*, U} \leq 0 \quad \forall z \in U_{ad}.\]
Unravelling the adjoint operator, we get the desired claim.
\end{proof}
An adjoint equation can be obtained, like in \cite{Kouri2018} and \cite{Kouri2019a}, by setting 
$(T_\tau'(u_\tau^*))^*\cJ'(y_\tau^*) = B^*p_\tau^* = {B(\cdot)^*}p_\tau^*(\cdot)$ 
with $p_\tau^*$ taking the role of the adjoint variable (this allows for the formulation of more explicit optimality conditions). Doing so, we get the following. 
\begin{prpstn}\label{prop:KKT-penalised}
{Let Assumption \ref{ass:onvarrho} hold.} There exists 
$(p_\tau^*, \pi_\tau^*) \in  L^{\tilde p}(\Omega;V) \times L^{p'}(\Omega)$ such that
\begin{subequations}\label{eq:tauSSOld}
\begin{align}
    {A(\omega)^*}p_\tau^*(\omega) + \frac 1\tau m_\tau'(y_\tau^*(\omega)-\psi(\omega))p_\tau^*(\omega) &= J_y(y_\tau^*(\omega),\omega) \quad \text{a.s.},\label{eq:ptau-old}\\
    \E[\langle B^* p_\tau^*, u_\tau^*-v \rangle_{U^*,U} \pi_\tau^*] + \varrho'(u_\tau^*)(u_\tau^*-v) &\leq 0 \quad \forall v \in U_{ad},\label{eq:tauSS3-old}\\
    \cR[g]-\cR[\cJ(y_\tau^*)] - \E[\pi_\tau^*(g-\cJ(y_\tau^*))] &\geq 0 \quad \forall g \in L^p(\Omega).\label{eq:tauSS4-old}
\end{align}
\end{subequations} 
\end{prpstn}
\begin{proof}
Define $p_\tau^*$ to satisfy \eqref{eq:ptau-old}. By Lemma \ref{lem:characterisationOfAdjointTtauPrime}, we have
\[(T'_{\tau,\omega}(u_\tau^*))^* J_y(y_\tau^* ,\omega)  = {B(\omega)^*}p_\tau^*(\omega).\]
Using this, the first term in inequality \eqref{eq:OC-penalized-as-VI} can be rewritten as
{
 \begin{align*}
    \E[\cJ'(y_\tau^*)(T'_\tau(u_\tau^*)(u_\tau^*-v))\pi_\tau^*] 
    &= \E[\langle \cJ'(y_\tau^*), T'_\tau(u_\tau^*)(u_\tau^*-v)\pi_\tau^*\rangle ]\\
    &= \E[\langle (T'_\tau(u_\tau^*))^*\cJ'(y_\tau^*), u_\tau^*-v\rangle_{U^*, U}\pi_\tau^* ]\\
    &= \E[\langle B^*p_\tau^*, u_\tau^*-v\rangle_{U^*, U}\pi_\tau^*, ]
 \end{align*} 
 }
 which gives \eqref{eq:tauSS3-old}.  The stated integrability on $p_\tau^*$ follows by Lemma \ref{lem:characterisationOfAdjointTtauPrime} from the regularity on its source term, see \eqref{ass:JyBoundedFinite}. Finally,  \eqref{eq:tauSS4-old} is simply  equivalent to the statement $\pi_\tau^* \in \partial\cR(\cJ(y_\tau^*))$.
\end{proof}

Now, in the sequel, we study the behavior of the corresponding sequences as $\tau \rightarrow 0.$ Inspecting \eqref{eq:tauSS3-old}, we see that limiting arguments will require a statement for the product $p_\tau^*\pi_\tau^*.$ It will turn out that our assumptions do not appear to provide strong convergence in either $p_\tau^*$ or $\pi_\tau^*$, making the identification of limits difficult. Therefore, we additionally pursue a different and slightly weaker formulation of the optimality conditions.  First, let us note that $q'$, the H\"older conjugate of $q$, satisfies\footnote{Indeed, we have $p' = \frac{p}{p-1}$ and $\tilde p = \frac{pq}{q-p}$ and so
\begin{align*}
    \frac{1}{p'} + \frac{1}{\tilde p} &= \frac{p-1}{p} + \frac{q-p}{pq} 
    = \frac{pq-p}{pq} = \frac{q-1}{q} 
    = \frac{1}{q'}.
\end{align*}}
\begin{equation}
\frac{1}{q'} =  \frac{1}{p'} + \frac{1}{\tilde p}.\label{eq:qprimeRelatedToPprimeEtc}
\end{equation}
and that $q' > 1$ (this is needed for applications of Banach--Alaoglu).

\begin{prpstn}\label{prop:KKT-penalised-2}
{Let Assumption \ref{ass:onvarrho} hold.} There exists $(q_\tau^*, \pi_\tau^*) \in L^{q'}(\Omega;V) \times L^{p'}(\Omega)$ such that
\begin{subequations}\label{eq:tauSSNew}
\begin{align}
    {A(\omega)^*}q_\tau^*(\omega) + \frac 1\tau m_\tau'(y_\tau^*(\omega)-\psi(\omega))q_\tau^*(\omega) &= J_y(y_\tau^*(\omega),\omega)\pi_\tau^*(\omega) \quad \text{a.s.},\label{eq:ptau}\\
    \E[\langle B^* q_\tau^*, u_\tau^*-v \rangle_{U^*,U}] + \varrho'(u_\tau^*)(u_\tau^*-v) &\leq 0 \quad \forall v \in U_{ad},\label{eq:tauSS3}\\
    \cR[g]-\cR[\cJ(y_\tau^*)] - \E[\pi_\tau^*(g-\cJ(y_\tau^*))] &\geq 0 \quad \forall g \in L^p(\Omega).\label{eq:tauSS4}
\end{align}
\end{subequations} 
\end{prpstn}
\begin{proof}
Define $q_\tau^*$ to satisfy \eqref{eq:ptau}. By Lemma \ref{lem:characterisationOfAdjointTtauPrime}, we have
\[(T'_{\tau,\omega}(u_\tau^*))^* J_y(y_\tau^* ,\omega) \pi_\tau^*(\omega) = {B(\omega)^*}q_\tau^*(\omega).\]
Using this, the first term in inequality \eqref{eq:OC-penalized-as-VI} can be rewritten as
 \begin{align*}
    \E[\cJ'(y_\tau^*)(T'_\tau(u_\tau^*)(u_\tau^*-v))\pi_\tau^*] 
    &= \E[\langle \cJ'(y_\tau^*), T'_\tau(u_\tau^*)(u_\tau^*-v)\rangle \pi_\tau^*]\\
    &= \E[\langle (T'_\tau(u_\tau^*))^*\cJ'(y_\tau^*)\pi_\tau^*, u_\tau^*-v\rangle_{U^*, U} ]\\
    &= \E[\langle B^* q_\tau^*, u_\tau^*-v\rangle_{U^*, U} ]
 \end{align*}
 which gives \eqref{eq:tauSS3}. Recalling \eqref{eq:qprimeRelatedToPprimeEtc}, we apply  H\"older's inequality to obtain
\[\left(\int_\Omega  |\pi_\tau^*|^{q'}\norm{J_y(y_\tau^*(\omega),\omega)}{V^*}^{q'} \,\mathrm{d}\pP(\omega) \right)^{\frac{1}{q'}} \leq C\norm{ \pi_\tau^*}{L^{p'}(\Omega)}\norm{\cJ'(y_\tau^*)}{L^{\tilde p}(\Omega;V^*)}, \]
which is finite. It follows that $q_\tau^*$ is also in this space.
\end{proof}

\begin{lmm}\label{lem:qtauEqualsPTauPiTau}
We have that in fact
\[q_\tau^* = p_\tau^*\pi_\tau^*.\]
\end{lmm}
\begin{proof}
If we multiply the equation \eqref{eq:ptau-old} for $p_\tau^*$ by the scalar $\pi_\tau^*$ and set $\hat q:= p_\tau^*\pi_\tau^*$, it immediately follows that
\[A(\omega)\hat q(\omega) + \frac 1\tau m_\tau'(y_\tau^*(\omega)-\psi(\omega))\hat q(\omega) = J_y(y_\tau^*(\omega),\omega)\pi_\tau^*(\omega).\]
Since $q_\tau^*$ also satisfies this equation, by uniqueness we must have that $q_\tau^* = \hat q.$
\end{proof}
\begin{prpstn}
If $p_\tau^*$ satisfies \eqref{eq:tauSSOld}, then the quantity $\hat q =p_\tau^*\pi_\tau^*$ satisfies \eqref{eq:tauSSNew}.
\end{prpstn}
This shows that the system \eqref{eq:tauSSOld} is in some sense stronger than \eqref{eq:tauSSNew}. More precisely, if we begin with \eqref{eq:tauSSOld}, we can obtain \eqref{eq:tauSSNew}. For the converse to hold (and thus for equivalence of the two systems) we would need $\pi_\tau^* \neq 0$ a.s., which is not the case in general.
\begin{proof}
If we define $\hat q:= p_\tau^*\pi_\tau^*$, we have, just like we argued above,
\[A(\omega)\hat q(\omega) + \frac 1\tau m_\tau'(y_\tau^*(\omega)-\psi(\omega))\hat q(\omega) = J_y(y_\tau^*(\omega),\omega)\pi_\tau^*(\omega)\]
is satisfied, i.e., we have \eqref{eq:ptau}, and \eqref{eq:tauSS3} is also immediate.
\end{proof}

    \section{Stationarity conditions}\label{sec:statConds}
MPECs may admit various types of stationarity conditions based on the assumptions and derivation techniques. In this section, we derive forms of weak and C-stationarity for the problem \eqref{eq:OCProblem}. Note that there are many refinements of the concept of C-stationarity and the terminology is  used somewhat inconsistently in the literature. Let us also remark that a number of related stationarity concepts can be derived for elliptic MPECs; see \cite{MR2515801, MR3796767} and also \cite{hintermuller2014several} for a systematic treatment of derivation techniques. Our approach is to derive strong stationarity conditions for the penalised problem \eqref{eq:penalisedOCProblem} and then to pass to the limit. First of all however, we must prove that \eqref{eq:penalisedOCProblem} is indeed a suitable penalisation for \eqref{eq:OCProblem}.

\subsection{Consistency of the approximation} 
We begin with an uniform estimate.
\begin{lmm}\label{lem:boundednessofyandutaus}{Under Assumption \ref{ass:onvarrho},} the following bound holds uniformly in $\tau >0$:
\begin{equation}
\label{eq:bound-y-u-tau}
    \norm{y_\tau^*}{L^{q}(\Omega;V)} + \norm{u^*_\tau}{U}\leq C.
\end{equation}
Thus, {there exists $(y^*, u^*) \in L^q(\Omega; V) \times U$ such that (for a subsequence that we do not distinguish)}
\begin{equation*}
    \begin{aligned}
    y_\tau^* &\weaklyto y^* &&\text{in $L^{q}(\Omega;V)$},\\
    u^*_\tau &\weaklyto u^* &&\text{in $U$}
    \end{aligned}
\end{equation*}
as $\tau \downarrow 0$.
\end{lmm}
\begin{proof}
The second bound in \eqref{eq:bound-y-u-tau} is due to Assumption \ref{ass:onvarrho} \ref{item:coerciveOrBddForTtauAndS}. 
Since $y_\tau^* = T_\tau(u_\tau^*)$, it follows by Proposition \ref{lem:boundsOnTtau} that $y_\tau^*$ is bounded uniformly in $L^{q}(\Omega;V)$. 
\end{proof}
In fact, we can do better for the state: since $y_\tau^* = T_\tau(u_\tau^*)$ and $u_\tau^* \weaklyto u^*$, we can immediately obtain strong convergence thanks to Proposition \ref{lem:fundamentalLemma}. 
\begin{lmm}\label{lem:yTauStrongConvergence}
We have {(for a subsequence that we do not distinguish)}
\[y_\tau^* \to y^* \qquad\text{in $L^{q}(\Omega;V)$}
 \]
 and $y^*$ solves in an a.s.~sense the variational inequality
 \begin{align*}
     y^*(\omega) \leq \psi(\omega), \qquad \langle A(\omega)y^*(\omega)-f(\omega)- B(\omega)u^*, y^*(\omega)-v \rangle &\leq 0 \qquad \forall v \in V : v \leq \psi(\omega).
\end{align*}
 \end{lmm}
{Having shown that $y^* = S(u^*)$, let us now show that $(y^*, u^*)$ indeed corresponds to an optimal pair of \eqref{eq:OCProblem}.}
\begin{prpstn}\label{prop:convergenceOfOpts}
{We have that $u^* \in U$ is a minimiser of \eqref{eq:OCProblem} with associated state $y^* \in L^q(\Omega;V).$}

\end{prpstn}
\begin{proof}
If $\bar u$ is an arbitrary minimiser of \eqref{eq:OCProblem} {(which exists by Proposition \ref{prop:existenceOC})} and $\bar y:= S(\bar u)$, then
\[\cR[\cJ(\bar y)] + \varrho(\bar u) \leq \cR[\cJ(S(u))] + \varrho(u) \quad \forall u \in U_{ad}.\]
In particular, for  $u=u^*$, we get 
\[\cR[\cJ(\bar y)] + \varrho(\bar u) \leq \cR[\cJ(y^*)] + \varrho(u^*).\]
On the other hand, we have for the penalised control problem
\[\cR[\cJ(y_\tau^*)] + \varrho(u_\tau^*) \leq \cR[\cJ(T_\tau(\hat u))] + \varrho(\hat u) \quad \forall \hat u \in U_{ad},\]
since $u_\tau^*$ is a minimiser of \eqref{eq:penalisedOCProblem}. In particular, with $\hat u$ selected as $\bar u$, we obtain
\[\cR[\cJ(y_\tau^*)] + \varrho(u_\tau^*) \leq \cR[\cJ(T_\tau(\bar u))] + \varrho(\bar u).\]
Using Proposition \ref{lem:fundamentalLemma}, we have $T_\tau(\bar u) \to S(\bar u)= \bar y$, and hence by continuity of {$\cJ$ (see \eqref{eq:JisCts})}, we have
\[\cR[\cJ(\bar y)] + \varrho(\bar u) \leq \cR[\cJ(y^*)] + \varrho(u^*) \leq \liminf_{\tau \to 0}\cR[\cJ(y_\tau^*)] + \varrho(u_\tau^*)  
\leq \cR[\cJ(\bar y)] + \varrho(\bar u),\]
where the first inequality is because $\bar u$ was assumed to be a minimiser and for the second inequality we used weak lower semicontinuity  of $\cR \circ \cJ \circ S$ (see the proof of Proposition \ref{prop:existenceOC}) and of $\varrho$ (see the discussion after Assumption \ref{ass:onvarrho}). We see that $\cR[\cJ(y^*)] + \varrho(u^*)$ coincides with the minimal value $\cR[\cJ(\bar y)] + \varrho(\bar u)$ and hence $u^*$ must be a minimiser.
\end{proof}

\subsection{Passage to the limit}

\begin{lmm}
\label{lem:bounded-dual-variables}
The following bound holds  uniformly in $\tau >0$ (for $\tau$ sufficiently small):
\begin{equation*}
    \norm{p_\tau^*}{L^{\tilde p}(\Omega;V)} + \norm{q_\tau^*}{L^{q'}(\Omega;V)} + \norm{\pi^*_\tau}{L^{p'}(\Omega)}\leq C.
\end{equation*}
Thus, we have (for a subsequence that we do not distinguish)
\begin{equation*}
    \begin{aligned}
    p_\tau^* &\weaklyto p^* &&\text{in $L^{\tilde p}(\Omega;V)$},\\
    q_\tau^* &\weaklyto q^* &&\text{in $L^{q'}(\Omega;V)$},\\
    \pi^*_\tau &\weaklyto \pi^* &&\text{in $L^{p'}(\Omega)$} \quad ( \weakstar  \text{ if } p'=\infty)
    \end{aligned}
\end{equation*}
as $\tau \downarrow 0.$
\end{lmm}
\begin{proof}
Test the equation \eqref{eq:ptau-old} with $p^*_\tau$ 
 and use the boundedness of $J_y$ in \eqref{ass:JyBoundedFinite} to obtain
\[C_a\norm{p_\tau^*(\omega)}{V} \leq \tilde{C}_1(\omega) + \tilde{C}_2\norm{y_\tau^*(\omega)}{V}^{q/\tilde{p}}.\]
By  Lemma \ref{lem:boundednessofyandutaus}, $y_\tau^*$ is bounded uniformly in $L^{q}(\Omega;V)$ and hence $p_\tau^*$ is bounded uniformly in $L^{\tilde{p}}(\Omega;V).$ 

For the bound on {$\pi_\tau^*$,} 
first observe that there exists a $\delta>0$ such that $\cR$ is Lipschitz with a Lipschitz constant $L$ on (the open ball) $B_\delta(\cJ(y^*)) \subset L^p(\Omega;V)$ since $\cR$ is locally Lipschitz.
The constant $L$ is obviously independent of $\tau$. For sufficiently small $\tau$, we get by continuity that $\cJ(y_\tau^*) \in B_{\delta}(\cJ(y^*))$.  It follows that $\cR$ is Lipschitz with the same Lipschitz constant $L$ near $\cJ(y_\tau^*)$ for $\tau$ small enough and hence, by \cite[Proposition 2.1.2 (a)]{Clarke} and the fact that $\pi_\tau^*  \in \partial\cR[\cJ(y_\tau^*)]$, we obtain the boundedness $\norm{\pi_\tau^*}{L^{p'}(\Omega)} \leq L$. 

The bound for $q_\tau^*$ follows similarly to the bound on $p_\tau^*$, using the fact that $\pi_\tau^*$ is bounded in $L^{p'}(\Omega)$.
\end{proof}

 \begin{rmrk}\label{rem:infty}
 If $q=\infty$, by the theory in \cite{Goldberg}, we would also need to have $p=\infty$ for the $C^1$ property for $\cJ$. In this case, $\partial \cR(y_\tau)$ is a subset of $L^\infty(\Omega)^*$, which is a space of measures (and $\pi_\tau^*$ would belong to this space). This is one reason why we postponed the case of $q=\infty$ to later work.
 \end{rmrk}

Unfortunately,  we cannot retrieve a (weak or strong) convergence pointwise a.s.~for $p_\tau^*$ by the same argument we used to prove Lemma \ref{lem:yTauStrongConvergence} because the limit may not be uniquely determined, as a (subsequence-dependent) multiplier would come into play.

\begin{rmrk}
 When $U_{ad}$ is the entire space, we obtain in the usual tracking-type non-stochastic setting that $p^*$ and the multipler  associated to the adjoint are both uniquely determined by $y^*$ and $u^*$. We explore this idea further.
 
 If $U_{ad} = U$, the inequality \eqref{eq:tauSS3-old} simplifies to
\[\left\langle \E[ \pi_\tau^*{B(\cdot)^*} p_\tau^*], v \right\rangle_{U^*,U}  + \varrho'(u_\tau^*)(v) = 0 \quad \forall v \in U.\]
Consider the case where $\varrho(u) := \frac{\nu}{2}\norm{u}{H}^2$ and $U = H$.  Then this further reduces to 
\[ \E[ \pi_\tau^*{B(\cdot)^*} p_\tau^*] + \nu u_\tau^* = 0\]
as an equality in $H^*$. Unfortunately, this does not seem to give us any strong convergence of $p_\tau^*$ or $\pi_\tau^*$ (unlike in the deterministic setting where the former would be available).
\end{rmrk}

Let us make an observation. If, {given $t \in (1,\infty)$}, $y \in L^t(\Omega;V)$ and $z \in L^{t'}(\Omega;V)$, by H\"older's inequality, we have
\[\left|\int_\Omega \langle A(\omega)y(\omega), z(\omega) \rangle \D \pP(\omega) \right| \leq C_b\norm{y}{L^t(\Omega;V)}\norm{z}{L^{t'}(\Omega;V)}\]
which implies that the associated operator considered between Bochner spaces, {which we denote by $\mathcal{A}$}, defined via the pairing
\begin{equation*}
   (\mathcal{A}y)(z) := \int_\Omega \langle A(\omega)y(\omega), z(\omega)\rangle \D \pP(\omega),
\end{equation*}
and mapping $L^t(\Omega;V)$ to $(L^{t'}(\Omega;V))^*$, is a bounded linear operator for any $t>1$.  Making free use of reflexivity and non-borderline exponent cases, we see that $\mathcal{A} \colon L^t(\Omega;V) \to L^t(\Omega;V^*)$ and the associated adjoint operator $\mathcal{A}^*\colon L^{t'}(\Omega;V) \to L^{t'}(\Omega;V^*)$ are both continuous.

For the ensuing analysis, it is useful to define {in an a.s sense}
\begin{align*}
\lambda_\tau^* (\omega):= \frac 1\tau m_\tau'(y_\tau^*(\omega)-\psi(\omega))p_\tau^*(\omega) &= J_y(y_\tau^*(\omega),\omega)-{A(\omega)^*}p_\tau^*(\omega),\\
\hat\lambda_\tau^* (\omega):= \frac 1\tau m_\tau'(y_\tau^*(\omega)-\psi(\omega))q_\tau^*(\omega) &= J_y(y_\tau^*(\omega),\omega)\pi_\tau^*-{A(\omega)^*}q_\tau^*(\omega),
\end{align*}
(these are  equalities in $V^*$). Note the relationship
\[\hat\lambda_\tau^* = \pi_\tau^*\lambda_\tau^*.\]

\begin{lmm}\label{lem:theLambdas}
{Under Assumption \ref{ass:onvarrho},} there exists $\lambda^* \in L^{\tilde p}(\Omega;V^*)$ and $\hat\lambda \in L^{q'}(\Omega;V^*)$ such that 
\begin{align*}
    \lambda^* &:= \cJ'(y^*) - A^*p^*, \qquad \hat \lambda^* := \cJ'(y^*)\pi^* - A^*q^*,
\end{align*}
and {(for a subsequence that we do not distinguish)}
\begin{align*}
\lambda_\tau^* &\weaklyto \lambda^* \quad\text{ in $L^{\tilde p}(\Omega;V^*)$}, \\ 
\hat\lambda_\tau^* &\weaklyto \hat\lambda^* \quad\text{ in $L^{q'}(\Omega;V^*)$}. 
\end{align*}
\end{lmm}
\begin{proof}
{In this proof we make use of Assumption \ref{ass:onvarrho} and its consequences.} Take $\varphi \in L^{\tilde{p}'}(\Omega;V)$ and consider
\begin{align*}
    \E[\langle \lambda_\tau^*, \varphi \rangle]  = \int_\Omega \langle J_y(y_\tau^*(\omega),\omega), \varphi(\omega) \rangle -  \langle A(\omega)\varphi(\omega), p_\tau^*(\omega) \rangle \D \pP(\omega).
\end{align*}
By  \eqref{eq:JyCtsGoodSpace}, we have $\mathcal{J}'(y_\tau^*) \to \mathcal{J}'(y^*)$ in $L^{\tilde p}(\Omega;V^*)$, taking care of the first term on the right-hand side above. The second follows trivially since we have $p_\tau^* \weaklyto p^*$ in $L^{\tilde p}(\Omega;V).$ 

For $\hat\lambda_\tau^*$, take now $\varphi \in L^q(\Omega;V)$ and consider
\begin{align*}
   \E[\langle \hat \lambda_\tau^*, \varphi \rangle] = \int_\Omega \pi_\tau^*(\omega) \langle J_y(y_\tau^*(\omega),\omega), \varphi(\omega) \rangle -  \langle A(\omega)\varphi(\omega), q_\tau^*(\omega) \rangle \D \pP(\omega).
\end{align*}
Since $\cJ$ is $C^1$, $\mathcal{J}' \colon L^q(\Omega;V) \to \mathcal{L}(L^q(\Omega;V), L^p(\Omega))$ is continuous and $\pi_\tau^* \weaklyto \pi^*$ in $L^{p'}(\Omega)$, thus we can pass to the limit in the first term on the right-hand side and using $q_\tau^* \weaklyto q^*$ in $L^{q'}(\Omega;V),$ we  can conclude.
\end{proof}

In preparation for the main result, recall  from \eqref{eq:qprimeRelatedToPprimeEtc} that
\[\frac{1}{q'} = \frac{1}{p'} + \frac{1}{\tilde p}.\]
Let us also define the \textit{inactive set} at $y^*$ by
\[\mathcal I^* := \{y^* < \psi\} := \{ (\omega, x) \in \Omega \times D: y^*(\omega)(x) < \psi(\omega)(x) \}.\]

\begin{thrm}[$\mathcal{E}$-almost weak stationarity]
\label{thm:ss-weaker}
Let Assumption \ref{ass:onvarrho} hold. For any local minimiser  {$u^* \in U_{ad}$ of \eqref{eq:OCProblem} with associated state $y^* \in  L^q(\Omega;V)$}, there exists  $\zeta^* \in L^q(\Omega;V)$, $q^* \in L^{q'}(\Omega;V)$,  $\pi^* \in L^{p'}(\Omega)$ and $\hat\lambda^* \in L^{q'}(\Omega;V^*)$ such that the following  system is satisfied: 
\begin{subequations}\label{eq:eaWeakStationarity}
\begin{align}
A(\omega)y^*(\omega) - f(\omega) - B(\omega)u^* +  \zeta^*(\omega) &= 0 \quad \text{a.s.},\label{eq:eqYZeta}\\
\zeta^*(\omega) \geq 0, \quad y^*(\omega) \leq \psi(\omega), \quad \langle \zeta^*(\omega), y^*(\omega)-\psi(\omega)\rangle &= 0 \quad \text{a.s.}\label{eq:ppQ},\\
    {A(\omega)^*}q^*(\omega) + \hat\lambda^*(\omega)  &= J_y(y^*(\omega),\omega)\pi^*(\omega) \quad \text{a.s.}, \label{eq:adjoint-smoothQ}\\
            \E[\langle B^* q^*, u^*-v \rangle_{U^*,U}] + \varrho'(u^*)(u^*-v) &\leq 0 \quad \forall v \in U_{ad},\label{eq:chiInequalityQ}\\ 
    \cR[g]-\cR[\cJ(y^*)] - \E[\pi^*(g-\cJ(y^*))] &\geq 0 \quad \forall g \in L^p(\Omega),\label{eq:piAndRiskMeasureNew}\\
 \E [\langle\zeta^*, q^*\rangle]  &= 0\quad \text{if $q=2$,} \label{eq:zetaQ}\\
        \limsup_{\tau \to 0}\E[\langle \hat \lambda_\tau^*, q_\tau^* \rangle] &\geq 0,\label{eq:limsupCond}\\
     \E[\langle \hat\lambda^*, y^*-\psi \rangle] &= 0,\\
       \forall \epsilon > 0, \exists E^\epsilon \subset \mathcal{I}^* \text{ with } |\mathcal{I}^* \setminus E^\epsilon| \leq \epsilon : \E[\langle \hat \lambda^*, v \rangle] &= 0 \quad \, \forall v \in L^{q}(\Omega;V)  :  v = 0 \text{ a.s.-a.e. on $\Omega\times D \setminus E^\epsilon$}.\label{eq:eaEaQ}
\end{align}
\end{subequations} 
\end{thrm}
\noindent Let us comment on this result.
\begin{enumerate}[label=(\roman*)]\itemsep0em 
\item The first two lines \eqref{eq:eqYZeta}--\eqref{eq:ppQ} encapsulate the well known complementarity form of the VI for $y^*$. 
\item We have not been able to obtain the sign condition
      $ \E[ \langle \hat\lambda^* , q^* \rangle] \geq 0$ 
(we only have \eqref{eq:limsupCond}). This is why we cannot call this a system of C-stationarity type; it resembles instead a weak stationarity system, only.  
\item In addition, we have \eqref{eq:zetaQ} only when $q=2$. We explain the complications that give rise to this (and a lack of further properties) after the proof of the theorem. 
\item Note here that Theorem \ref{thm:ss-weaker} holds true for local minimisers whereas Proposition \ref{prop:convergenceOfOpts} argues for global minimisers of the the respective problems.
\end{enumerate}
The theorem essentially follows by passing to the limit in the $q_\tau^*$ system \eqref{eq:tauSSNew}. Before we get to that, it is useful and instructive to first pass to the limit in the adjoint equation in the $p_\tau^*$ system \eqref{eq:tauSSOld} and to derive properties of its associated quantities.
\begin{prpstn}[Weak $\mathcal{E}$-almost C-stationarity]\label{thm:ss}
Let Assumption \ref{ass:onvarrho} and
\begin{equation}
   q \leq 2p \quad\text{or} \quad p \geq 2\label{ass:furtherRestrictionOnPq}
\end{equation}
hold. 
For any local minimiser {$u^* \in U_{ad}$ of \eqref{eq:OCProblem} with associated state $y^* \in L^q(\Omega;V)$}, there exists $\zeta^* \in L^q(\Omega;V^*)$, 
$p^* \in L^{\tilde{p}}(\Omega;V)$  and $\lambda^* \in L^{\tilde p}(\Omega;V^*)$ 
such that the following  system is satisfied:
\begin{subequations}\label{eq:ecsSytem}
\begin{align}
A(\omega)y^*(\omega) - f(\omega) - B(\omega)u^* +  \zeta^*(\omega) &= 0 \quad \text{a.s.},\label{eq:eqYZeta-new}\\
\zeta^*(\omega) \geq 0, \quad y^*(\omega) \leq \psi(\omega), \quad \langle \zeta^*(\omega), y^*(\omega)-\psi(\omega)\rangle &= 0 \quad \text{a.s.},\label{eq:pp-new}\\
    {A(\omega)^*}p^*(\omega) + \lambda^*(\omega)  &= J_y(y^*(\omega),\omega) \quad \text{a.s.},\label{eq:adjoint-smooth}\\
        \E [\langle\zeta^*, p^*\rangle] &= 0,\label{eq:zetaP}\\
       \E[ \langle \lambda^* , p^* \rangle] &\geq 0,\label{eq:TT}\\
     \E[\langle \lambda^*, y^*-\psi \rangle] &= 0,\label{eq:complementarityLambda}\\
       \forall \epsilon > 0, \exists E^\epsilon \subset \mathcal{I}^* \text{ with } |\mathcal{I}^* \setminus E^\epsilon| \leq \epsilon : \E[\langle \lambda^*, v \rangle] &= 0 \, \forall v \in L^{\tilde p'}(\Omega;V)  : v = 0 \text{ a.s.-a.e. on $\Omega\times D \setminus E^\epsilon$}.\label{eq:eaEa}
\end{align}
\end{subequations}
\end{prpstn}
\noindent Before we prove this result, some remarks are in order.
\begin{enumerate}[label=(\roman*)]\itemsep0em 
    \item The system \eqref{eq:ecsSytem}
    resembles the so-called $\mathcal{E}$-almost C-stationarity system in the deterministic setting \cite{MR2515801}. However, in contrast to what we expected from the deterministic setting, we are unable to show that
\[\E[ \langle \zeta^*, (p^*)^+ \rangle] = \E[ \langle \zeta^*, (p^*)^- \rangle] = 0,\]
see Remark \ref{rem:unabletoshow}. This is why we added the adjective `weak' to refer to the system.
\item We have not been able to relate the adjoint $p^*$ with the control $u^*$ nor $\pi^*$. This is why we have presented the system \eqref{eq:eaWeakStationarity} as our main result where such a relationship between the adjoint and control is available. 
    \item Taking a cue from Lemma \ref{lem:qtauEqualsPTauPiTau}, we would ideally like to identify as $q^*$ with $p^*\pi^*$; this is a major issue. We show in Proposition \ref{prop:identification}  that this identification does hold in certain circumstances.
\item Regarding relations such as \eqref{eq:TT}, we cannot say anything about the duality products in an a.s. sense (i.e., without the expectation present) because we do not have the desired convergences of elements such as $\lambda_\tau^*(\omega)$ and $p_\tau^*(\omega)$ for fixed $\omega$.
\item If \eqref{ass:furtherRestrictionOnPq} is not available, we do not get \eqref{eq:zetaP}--\eqref{eq:TT} (but the remaining statements still hold). This is essentially because we need the estimate \eqref{eq:lambdaPBound} and we have it under \eqref{ass:furtherRestrictionOnPq} as it implies that Lemma \ref{lem:tildePGeqTwo} is applicable, which yields $\tilde p \geq 2$.
\end{enumerate}
\begin{proof}[Proof of Proposition \ref{thm:ss}]
The proof is in part similar to that of \cite[Theorem 3.4]{MR2822818} but more delicate due to the low regularity convergence results.

  \begin{enumerate}[wide, labelwidth=!, labelindent=0pt]
\item \emph{The adjoint equation \eqref{eq:adjoint-smooth}.} This is a simple consequence of Lemma \ref{lem:theLambdas}. 
\item\label{step:lambdaP} \emph{Sign condition \eqref{eq:TT} on the product.} 
Observe that
\begin{equation}
\E[\langle  \lambda_\tau^*, p_\tau^* \rangle] \leq \norm{\lambda_\tau^*}{L^2(\Omega;V^*)}\norm{p_\tau^*}{L^2(\Omega;V)} \leq C\norm{\lambda_\tau^*}{L^{\tilde p}(\Omega;V^*)}\norm{p_\tau^*}{L^{\tilde p}(\Omega;V)}\label{eq:lambdaPBound}
\end{equation}
since (as remarked above) $L^{\tilde p}(\Omega;V) \cts L^2(\Omega;V)$; so the left-hand side is well defined. Now, since $\langle \lambda_\tau^*, p_\tau^* \rangle \geq 0$ a.s.~(using the definition), we have
\begin{align*}
    0 &\leq \limsup_{\tau \to 0} \E[\langle \lambda_\tau^*, p_\tau^* \rangle]\\
    &\leq \limsup_{\tau \to 0} \E[ \mathcal{J}'(y_\tau^*)( p_\tau^*)] - \liminf_{\tau \to 0} \E[\langle A^*p_\tau^*, p_\tau^* \rangle]\\
    &\leq \E[\langle \mathcal{J}'(y^*), p^* \rangle] - \E[\langle A^*p^*, p^* \rangle]\\
    &= \E[\langle \lambda^*, p^* \rangle].
\end{align*}
Here, to derive the penultimate line, we used the strong convergence of $y_\tau^*$, the continuity of $\mathcal{J}'$ into   $L^{\tilde p}(\Omega;V^*)$ from \eqref{eq:JyCtsGoodSpace} in combination with the fact that $\tilde p \geq \tilde p'$ (see  Lemma \ref{lem:trivialLemma}), which implies that $p_\tau^* \weaklyto p^*$ in $L^{\tilde p'}(\Omega;V)$ as well as the weak lower semicontinuity 
of the bounded, coercive bilinear form $\E[\langle A^* \cdot, \cdot)\rangle]\colon L^2(\Omega;V) \times L^2(\Omega;V) \to \mathbb{R}$. {The expectation in the bilinear form is finite by Lemma \ref{lem:tildePGeqTwo} as we argued above.}

\item \emph{Relation \eqref{eq:complementarityLambda} between multiplier and $y^*-\psi$.}  We have
\begin{align*}
    \langle \lambda_\tau^*(\omega), (y_\tau^*(\omega)-\psi(\omega))^- \rangle = \frac 1\tau \langle m_\tau'(y_\tau^*(\omega)-\psi(\omega))p_\tau^*(\omega), (y_\tau^*(\omega)-\psi(\omega))^- \rangle = 0
\end{align*}
because the duality product above is just the integral over the domain and $m_\tau'$ vanishes on the negative line. Taking the expectation 
and using the strong convergence $y_\tau^* \to y^*$ in $L^q(\Omega;V)$ {from Lemma \ref{lem:yTauStrongConvergence}}, which certainly implies strong convergence of its positive (and negative) part in $L^q(\Omega;V)$, and using also $\lambda_\tau^* \weaklyto \lambda^*$ in $L^{q'}(\Omega;V^*)$ by \eqref{eq:tildePGreaterThanP}, we can pass to the limit and then realising that $(y^*-\psi)^- = -(y^*-\psi)$, we find the desired condition.

\item \emph{$\mathcal{E}$-almost statement \eqref{eq:eaEa}.} Since $y_\tau^* \to y^*$ in $L^q(\Omega; L^2(D))$,  due to the identification $L^1(\Omega; L^1(D)) \cong L^1(\Omega \times D)$, we have $y_\tau^*-\psi \to y^*-\psi$ pointwise a.s.-a.e.~in $\Omega \times D =: \Omega_D$ for a subsequence that we do not distinguish. Let the measure of $\mathcal{I}^*$ be positive; otherwise nothing needs to be shown. Then take $z \in  \Omega_D$ such that $y^*(z)-\psi(z) < 0$, then there exists a $\hat \tau=\hat \tau(z)$ such that if $\tau \leq \hat \tau$, then
\[y_\tau^*(z)-\psi(z) \leq \frac12 (y^*(z)-\psi(z)) < 0\]
and hence $\tau^{-1}m_\tau'(y_\tau^*(z)-\psi(z)) = 0$ for $\tau \leq \hat \tau$. That is, $\tau^{-1}m_\tau'(y_\tau^*(z)-\psi(z)) \to 0$ pointwise a.s.-a.e.~on $\{y^* < \psi\}$ and by Egorov's theorem, for every $\epsilon > 0$, there exists $B^\epsilon \subset \{y^* < \psi\}$ with $|B^\epsilon| < \epsilon$ such that this convergence also holds uniformly on $\{y^* < \psi\}\setminus B^\epsilon$. 

Take $v \in L^{\tilde p'}(\Omega;V)$ with $v=0$ a.s.-a.e.~on $\{y^*=\psi\}\cup B^\epsilon \subset \Omega_D$. By the uniform convergence, for any $\gamma > 0$, there exists $\bar\tau$ such that if $\tau \leq \bar\tau$,
\begin{align*}
\left|\int_\Omega \langle \lambda_\tau^*, v \rangle \D \pP(\omega)\right| = \left|\int_{\{y^* < \psi \} \cap (B^\epsilon)^c} \frac 1\tau  m_\tau'(y_\tau-\psi)p_\tau^*v  \D \pP(\omega) \right| \leq \gamma\norm{p_\tau^* v}{L^1(\Omega_D)} &= \gamma \norm{p_\tau^*v}{L^1(\Omega;L^1(D))}.
\end{align*}
The norm on the right-hand side is bounded independently of $\tau$ and the left-hand side converges to $|\E[ \langle \lambda^*, v \rangle ]|$ (thanks to $\lambda^*_\tau \weaklyto \lambda^*$ in $L^{\tilde p}(\Omega;V^*)$), thus giving 
\begin{align*}
\left|\E[\langle \lambda^*, v \rangle ] \right| \leq  C\gamma
\end{align*}
for a constant $C>0$. Since this holds for every $\gamma$, we obtain \eqref{eq:eaEa} (simply set $E^\epsilon := \mathcal{I}^*\setminus B^\epsilon$).

\item\label{step:zetaP} \emph{Relation \eqref{eq:zetaP} between $\zeta^*$ and $p^*$.} 
Define 
\[\zeta_\tau^* := \frac 1\tau m_\tau(y_\tau^*-\psi) = f+Bu_\tau^* - Ay_\tau^*\]
which satisfies $\zeta_\tau^* \to \zeta^*$ in $L^q(\Omega;V^*)$. Hence, since $q \geq 2$ we have
\[\E[\langle \zeta_\tau^*, y_\tau^*-\psi \rangle] \to \E[\langle \zeta^*, {y^*}-\psi \rangle] = 0,\]
with the equality due to \eqref{eq:ppQ}. Recall the definition of $m_\tau$ from \eqref{eq:mrhoHK}. Introducing 
{the two sets
\begin{align*}
M_1(\tau) &:= \{ (\omega, x) \in \Omega \times D : 0 \leq y_\tau^*(\omega, x) - \psi(\omega, x) < \tau\},\\
M_2(\tau) &:= \{(\omega, x) \in \Omega \times D : y_\tau^*(\omega, x) - \psi(\omega, x) \geq \tau\},
\end{align*}
} by the convergence above, we find 
\begin{align}
\nonumber \E[ \langle \zeta_\tau^*, y_\tau^* - \psi\rangle]&=\frac 1\tau \int_\Omega \int_D m_\tau(y_\tau^*-\psi)(y_\tau^* - \psi) \D x \D \pP(\omega)\\
&= \frac 1\tau  \int_{M_1(\tau)}\frac{(y_\tau^* - \psi)^3}{2\tau}\D x \D \pP(\omega) + \frac 1\tau \int_{M_2(\tau)}\left(y_\tau^* - \psi-\frac \tau 2\right)(y_\tau^* - \psi)\D x \D \pP(\omega)\label{eq:second}\\
\nonumber &\to 0,
\end{align}
and as both integrands in \eqref{eq:second} are non-negative, each integral must individually converge to zero too. Hence, using $L^2(\Omega;L^2(D)) \cong L^2(\Omega \times D)$,
\begin{equation}\label{eq:pr1}
\norm{\frac{\chi_{M_1(\tau)}(y_\tau^* - \psi)^{\frac 32}}{ \tau  }}{L^2(\Omega;L^2(D))} \to 0\qquad\text{and}\qquad \norm{\frac{\chi_{M_2(\tau)}(y_\tau^* - \psi-\frac\tau 2)}{\sqrt{\tau}}}{L^2(\Omega;L^2(D))} \to 0,
\end{equation}
where for the second convergence we used the fact that $y_\tau^* - \psi \geq y_\tau^* - \psi-\tau\slash 2 \geq 0$. 
We calculate
\begin{align}
\nonumber \E[ \langle \zeta_\tau^*, p_\tau^* \rangle] &= \frac 1\tau \int_\Omega  \int_{M_1(\tau)}\frac{(y_\tau^*-\psi)^2}{2\tau}p_\tau^*  + \frac 1\tau  \int_\Omega \int_{M_2(\tau)}\left(y_\tau^*-\psi-\frac \tau 2\right)p_\tau^*\\
\nonumber &=  \frac 12 \int_\Omega \int_D \chi_{M_1(\tau)}\frac{(y_\tau^*-\psi)^{3\slash 2}}{\tau}\frac{(y_\tau^*-\psi)^{1\slash 2}}{\tau}\chi_{M_1(\tau)}p_\tau^*  +  \int_\Omega \int_D \frac{\chi_{M_2(\tau)}\left(y_\tau^*-\psi-\frac \tau 2\right)}{\sqrt{\tau}}\frac{\chi_{M_2(\tau)}p_\tau^*}{\sqrt{\tau}}\\
&=  \frac 12 \left(\chi_{M_1(\tau)}\frac{(y_\tau^*-\psi)^{3\slash 2}}{\tau}, \frac{(y_\tau^*-\psi)^{1\slash 2}}{\tau}\chi_{M_1(\tau)}p_\tau^*\right)   +  \left(\frac{\chi_{M_2(\tau)}\left(y_\tau^*-\psi-\frac \tau 2\right)}{\sqrt{\tau}}, \frac{\chi_{M_2(\tau)}p_\tau^*}{\sqrt{\tau}}\right)\label{eq:productXiP},
\end{align}
where the inner products in the final line are in $L^2(\Omega;L^2(D))$. Now, using \eqref{eq:pr1}, the first term in each inner product above converges to zero and hence the above right-hand side will converge to zero if we are able to show that the second term in each inner product remains bounded.

From \eqref{eq:lambdaPBound}, making use of the boundedness of $\lambda_\tau^*$ and $p_\tau^*$, we derive 
\begin{align*}
\nonumber C &\geq |\E[\langle \lambda_\tau^*, p_\tau^* \rangle]|\\
\nonumber  &=\frac 1\tau \left|\int_\Omega\int_D  m_\tau'(y_\tau^*-\psi)(p_\tau^*)^2 \D x \D \pP(\omega) \right|\\
&=\frac{1}{\tau}\int_\Omega\int_{D}\chi_{M_1(\tau)}\frac{y_\tau^*-\psi}{\tau}(p_\tau^*)^2  \D x \D \pP(\omega) + \frac{1}{\tau}\int_\Omega\int_{D}\chi_{M_2(\tau)}(p_\tau^*)^2  \D x \D \pP(\omega).
\end{align*}
Both of the terms on the right-hand side are individually bounded uniformly in $\tau$ as the integrands are non-negative. This fact then implies from \eqref{eq:productXiP} that 
\[\E[ \langle\zeta^*, p^*\rangle]= 0.\qedhere\]
\end{enumerate}
\end{proof}

\begin{rmrk}\label{rem:unabletoshow}
Replacing $p_\tau^*$ by $(p_\tau^*)^+$ in \eqref{eq:productXiP} and in the above calculation, we obtain in the same way as in the proof above
\[\E[ \langle \zeta_\tau^*, (p_\tau^*)^+ \rangle] = 0.\]
From here, since we have only weak convergence of $p_\tau^*$, we cannot say that $(p_\tau^*)^+$ converges (weakly) to $(p^*)^+$ and pass to (and be able to identify) the limit in the above. We can only deduce 
\[\lim_{\tau \to 0}\E[ \langle \zeta_\tau^*, (p_\tau^*)^+ \rangle] = 0.\]
\end{rmrk}

We now prove the main result of this paper.
\begin{proof}[Proof of Theorem \ref{thm:ss-weaker}]
We can in part capitalise on the proof of Proposition \ref{thm:ss} but first we begin with the VI for $q^*$.

  \begin{enumerate}[wide, labelwidth=!, labelindent=0pt]
\item \emph{The VI \eqref{eq:chiInequalityQ} relating the control to a multiplier.} Since the first term in the inequality \eqref{eq:tauSS3} contains a product of a strongly convergent sequence with another product of two weakly convergent sequences, we need to use some compactness to pass to the limit here.  First, we write the first term of that inequality as
\begin{align*}
    \E[\langle B^* q_\tau^*, u_\tau^*-v \rangle_{U^*,U} ] 
    &= \int_\Omega \langle q_\tau^*(\omega), B(\omega)(u_\tau^*-v) \rangle \;\mathrm{d}\pP(\omega).
\end{align*}
Now, by \eqref{eq:DCTConsequence},  we have $B(u_\tau^*-v) \to B(u^*-v)$ 
in $L^{q}(\Omega;V^*)$. This  and the weak lower semicontinuity of Assumption \ref{ass:onvarrho} \ref{item:wlscOfVarrhoPrime} on $\varrho'$ allows us to pass to the limit in the inequality \eqref{eq:tauSS3}.

\item \emph{Inequality \eqref{eq:piAndRiskMeasureNew}  for the risk measure.} 
It is easy to pass to the limit in the inequality \eqref{eq:tauSS4} since $\pi_\tau^*$ converges weakly in $L^{p'}(\Omega)$ and by continuity (see \eqref{eq:JisCts}), $\mathcal{J}(y_\tau^*)$ converges strongly in $L^p(\Omega)$; thus we get the inequality after making use of the lower semicontinuity of $\cR$.

\item \emph{The statements \eqref{eq:adjoint-smoothQ}, \eqref{eq:zetaQ}--\eqref{eq:eaEaQ}.} The proof of the remaining statements in \eqref{eq:eaWeakStationarity} is more or less identical to the proof of Proposition \ref{thm:ss}. Let us point out the changes.  
Since $\hat\lambda_\tau^*$ and $q_\tau^*$ are both uniformly bounded in $L^{q'}$ with respect to $\Omega$, if $q' \geq 2$, we would have that the left-hand side of \eqref{eq:lambdaPBound} (with $q_\tau^*$ instead) is bounded uniformly. However, recalling that we assumed in \eqref{eq:defnOfq} that $q \geq 2$, we must choose $q=q'=2$ to avail of the estimate. In this case 
 \eqref{eq:lambdaPBound} (with $p_\tau^*$ replaced by $q_\tau^*$) still holds, and the right-hand side is bounded uniformly:
\[\E[\langle  \hat \lambda_\tau^*, q_\tau^* \rangle] \leq \lVert{\hat \lambda_\tau^*}\rVert_{L^2(\Omega;V^*)}\norm{q_\tau^*}{L^2(\Omega;V)} \leq C.\]
Thus, step \ref{step:zetaP} of the proof of Proposition \ref{thm:ss} is still valid and we obtain \eqref{eq:zetaQ}.

\item \emph{Conclusion.} 
We are left to show that the stationarity system in fact holds for all local minimisers, not just for a cluster point of {$\{u_\tau^*\}$}. The argument is classical. Suppose that {$\hat u$ is an arbitrary local minimiser with associated state $\hat y$}, so there exists a ball $B^U_{\gamma}(\hat u)$ in $U$ of radius $\gamma$ on which it is the minimiser. We modify \eqref{eq:penalisedOCProblem} as follows: 
\begin{equation}
     \min_{u \in U_{ad} \cap B^U_{\gamma}(\hat u)} \cR[\cJ(T_\tau(u))] +\varrho(u) + \norm{u-\hat u}{U}^2\label{eq:ocProblemPenBarbu},
\end{equation}
and we denote by {$\bar u_\tau$ a minimiser of this problem with $\bar y_\tau = T_\tau(\bar y_\tau)$}. Let us denote the non-reduced functional appearing above as
\[{\hat F(y,u)} := \cR[\cJ(y)] +\varrho(u) + \norm{u-\hat u}{U}^2.\]
From ${\hat F(\bar y_\tau,  \bar u_\tau)} \leq {\hat F(T_\tau(\hat u), \hat u)}$ and $T_\tau(\hat u) \to S(\hat u) = \hat y$ (recall Proposition \ref{lem:fundamentalLemma}), we have 
\begin{equation}
\limsup_{\tau \to 0}{\hat F(\bar y_\tau, \bar u_\tau)} \leq \cR[\cJ(\hat y)] +\varrho(\hat u).\label{eq:forref1}
\end{equation}
On the other hand, due to Assumption \ref{ass:onvarrho} \ref{item:coerciveOrBddForTtauAndS}, we obtain the existence of $v \in U$ such that (for a subsequence that we will not distinguish) $\bar u_\tau \weaklyto v$ in $U$ and $\bar y_\tau \to S(v) =: z$ in $L^q(\Omega;V)$, giving (by the identity $\limsup(a_n) + \liminf(b_n) \leq \limsup(a_n + b_n)$ and using weak lower semicontinuity) 
\begin{equation}
\limsup_{\tau \to 0}{\hat F(\bar y_\tau, \bar u_\tau)} \geq \cR[\cJ(z)] +\varrho(v)+ \limsup_{\tau \to 0}\norm{\bar u_\tau-\hat u}{U}^2 \geq \cR[\cJ(\hat y)] +\varrho(\hat u) + \limsup_{\tau \to 0}\norm{\bar u_\tau-\hat u}{U}^2, \label{eq:forref2}    
\end{equation}
with the last inequality because {$\hat u$} is a local minimiser and $v$ remains in $B_\gamma^U(\hat u)$. Combining {the two inequalities \eqref{eq:forref1} and \eqref{eq:forref2}} shows that $\hat u = v$ and $\bar u_\tau \to \hat u$ in $U$. 

{Now take an arbitrary $w \in B_{\gamma\slash 2}^U(\bar u_\tau)$. Since $\bar u_\tau \to \hat u$, it follows that if $\tau$ is sufficiently small, say $\tau \leq \tau_0$, $\bar u_\tau \in B_{\gamma \slash 2}^U(\hat u)$. We take $\tau \leq \tau_0$ from now on. Using the triangle inequality, 
\begin{align*}
    \norm{w- \hat u}{U} &\leq \norm{w- \bar u_\tau}{U} + \norm{\bar u_\tau - \hat u}{U}\\
    &\leq \gamma,
\end{align*}
i.e., $w \in B_{\gamma}^U(\hat u)$, and hence we deduce (invoking the fact that $\bar u_\tau$ is a minimiser on $U_{ad} \cap B_\gamma^U(\hat u)$)
\begin{align*}
    \hat F(\bar y_\tau, \bar u_\tau) = \hat F(T_\tau(\bar u_\tau), \bar u_\tau) \leq \hat F(T_\tau(w), w) \quad \forall w \in B_{\gamma\slash 2}^U(\bar u_\tau).
\end{align*}
This statement proves that $\bar u_\tau$ is a local minimiser of the problem
\[     \min_{u \in U_{ad}}  \cR[\cJ(T_\tau(u))] +\varrho(u) + \norm{u-\hat u}{U}^2.\]
The same arguments as above can be used to derive stationarity conditions for this problem} and in passing to the limit in those conditions, we will find that $(\hat y, \hat u)$ satisfies the same conditions as above.\qedhere
\end{enumerate}
\end{proof}
Inspecting the proof, we see that even if $q=2,$ we cannot pass to the limit in $\E[\langle \hat \lambda_\tau^*, q_\tau^* \rangle]$ (as in step \ref{step:lambdaP} of the proof of Proposition \ref{thm:ss}) due to the quantity $\E[\pi_\tau^* \cJ'(y_\tau^*)(q_\tau^*)]$, where we again have a product of weakly convergent sequences. If $q>2$, then $\E[\langle \hat \lambda_\tau^*, q_\tau^*\rangle]$ may not be finite since it is the integral of two elements that are only known to be $q'$-integrable with respect to $\Omega$ and $q'  \in (1,2)$. Therefore, a version of the estimate \eqref{eq:lambdaPBound} may not even exist, which, at least via our method of proof, rules out \eqref{eq:zetaQ}.

\begin{prpstn}\label{prop:identification}
Suppose $\partial \cR(z) = \{\cR'(z)\}$ holds for all $z \in L^p(\Omega)$ and $\cR'$ satisfies the property
\begin{equation}
    z_n \to z \text{ in $L^p(\Omega;V)$ } \implies \cR'(z_n) \to \cR'(z) \text{ in $L^{p'}(\Omega)$}.\label{ass:onRBeingNice}
\end{equation}
Then
\[\pi_\tau^* \to \pi^* \text{ in $L^{p'}(\Omega)$},\]
and
we can identify
\[q^* = \pi^*p^*.\]
Hence \eqref{eq:chiInequalityQ} can be strengthened to
\begin{align*}
        \E[\langle B^* p^*, u^*-v \rangle_{U^*,U}\pi^*] + \varrho'(u^*)(u^*-v) &\leq 0 \quad \forall v \in U_{ad}.
\end{align*}
\end{prpstn}
If $\cR$ is G\^ateaux differentiable, then the subdifferential reduces to a singleton \cite[Proposition 2.3.6 (d)]{Clarke} as required above. The assumption \eqref{ass:onRBeingNice} is a strong one, but if for example $\cR$ is continuously Fr\'echet differentiable, then it holds. Clearly, the case $\cR := \E$ meets all of these conditions. {For nontrivial examples, the results in \cite{Kouri2019} demonstrate how to construct convex risk measures from, e.g., coherent risk measures, using epiregularization. Depending on the type of smoothing functional, these new measures are continuously Fr\'echet differentiable.}
\begin{proof}
We now have that $\pi_\tau^* = \cR'(\cJ(y_\tau^*))$. Since $y_\tau^* \to y^*$, using {assumption \eqref{ass:onRBeingNice} and the continuity of $\cJ$}, we obtain the strong convergence $\pi_\tau^* \to \pi^*$ in $L^{p'}(\Omega)$.

Take $\varphi \in L^\infty(\Omega;V^*)$ and suppose that $p'< \infty$. We have then
\begin{align*}
    \norm{\pi_\tau^*\varphi - \pi^*\varphi}{L^{p'}(\Omega;V^*)}^{p'} 
    = \int_\Omega |\pi_\tau^*-\pi^*|^{p'}\norm{\varphi}{V^*}^{p'} \D \pP(\omega)
    \leq \norm{\varphi}{L^\infty(\Omega;V^*)}^{p'}\norm{\pi_\tau^*-\pi^*}{L^{p'}(\Omega)}^{p'}     &\to 0.
\end{align*}
Now, using \eqref{eq:tildePGreaterThanP}, we have $p_\tau^* \weaklyto p^*$ in $L^p(\Omega;V)$ and thus 
\begin{align*}
    \E[\langle p_\tau^*\pi_\tau^*, \varphi \rangle] &= \E[ \langle p_\tau^*, \pi_\tau^*\varphi \rangle] \to  \E[ \langle p^*, \pi^*\varphi \rangle ]=  \E[ \langle p^*\pi^*, \varphi \rangle].
\end{align*}
This shows that $p_\tau^*\pi_\tau^* \weakstar p^*\pi^*$ in $(L^\infty(\Omega;V^*))^*$ (and also in some space of vector-valued distributions). But since we already know that $p_\tau^*\pi_\tau^* \weaklyto q^*$ in $L^{q'}(\Omega;V)$, we must have $q^* = p^*\pi^*.$ 

If $p'=\infty$, obvious modifications to the above yield the same conclusions.
\end{proof}
It is well worth stating the full system that we obtain under the conditions of the above proposition. As mentioned, we do obtain this when $\cR$ is chosen to be $\E$, i.e., the risk-neutral case is covered.
\begin{crllr} \label{cor:fullSystem}
Let Assumption \ref{ass:onvarrho}, \eqref{ass:furtherRestrictionOnPq}, and the assumptions of Proposition \ref{prop:identification} hold. For any local minimiser {$u^* \in U_{ad}$  of \eqref{eq:OCProblem} with associated state $y^* \in  L^q(\Omega;V)$}, there exists  $\zeta^* \in L^q(\Omega;V)$, $p^* \in L^{\tilde p}(\Omega;V)$,  $\pi^* \in L^{p'}(\Omega)$ and $\lambda^* \in L^{\tilde p}(\Omega;V^*)$ such that the following  system is satisfied:
\begin{subequations} 
\begin{align}
A(\omega)y^*(\omega) - f(\omega) - B(\omega)u^* +  \zeta^*(\omega) &= 0 \quad \text{a.s.},\\
\zeta^*(\omega) \geq 0, \quad y^*(\omega) \leq \psi(\omega), \quad \langle \zeta^*(\omega), y^*(\omega)-\psi(\omega)\rangle &= 0 \quad \text{a.s.},\\
  {A(\omega)^*}p^*(\omega) + \lambda^*(\omega)  &= J_y(y^*(\omega),\omega) \quad \text{a.s.},\\
              \E[\langle B^* p^*, u^*-v \rangle_{U^*,U}\pi^*] + \varrho'(u^*)(u^*-v) &\leq 0 \quad \forall v \in U_{ad},\\ 
                  \cR[g]-\cR[\cJ(y^*)] - \E[\pi^*(g-\cJ(y^*))] &\geq 0 \quad \forall g \in L^p(\Omega),\\
        \E [\langle\zeta^*, p^*\rangle] &= 0,\\
             \E [\pi^*\langle\zeta^*, p^*\rangle]  &= 0\quad\text{if $q=2$},\tag{\theequation'}\\
       \E[ \langle \lambda^* , p^* \rangle] &\geq 0,\\
     \E[\langle \lambda^*, y^*-\psi \rangle] &= 0,\\
          \E[\pi^*\langle \lambda^*, y^*-\psi \rangle] &= 0,\tag{\theequation'}\\ 
       \forall \epsilon > 0, \exists E^\epsilon \subset \mathcal{I}^* \text{ with } |\mathcal{I}^* \setminus E^\epsilon| \leq \epsilon : \E[\langle \lambda^*, v \rangle] &= 0 \, \forall v \in L^{\tilde p'}(\Omega;V)  : v = 0 \text{ a.s.-a.e. on $\Omega\times D \setminus E^\epsilon$},\\
       \forall \epsilon > 0, \exists E^\epsilon \subset \mathcal{I}^* \text{ with } |\mathcal{I}^* \setminus E^\epsilon| \leq \epsilon : \E[\pi^*\langle \lambda^*, v \rangle] &= 0 \, \forall v \in L^{q}(\Omega;V)  : v = 0 \text{ a.s.-a.e. on $\Omega\times D \setminus E^\epsilon$}.\tag{\theequation'}
\end{align}
\end{subequations}
\end{crllr}
\section{Numerical example}
\label{sec:numerics}
In this section, we take a specific example for which lack of strict complementarity holds; this gives rise to a genuinely nonsmooth solution map $S$ for the underlying VI. As a proof of concept, we use a stochastic approximation algorithm (rather than developing an algorithm tailored to this specific problem class). 

\subsection{Problem formulation}
For the numerical experiments, we focus on a particular realisation of  problem \eqref{eq:OCProblem} subject to the random VI \eqref{eq:VIProblemFormal}, namely a modification of the example in Section \ref{sec:example}. We use $D=(0,1)\times (0,1)$, $U_{ad}=L^2(D)$, the tracking-type function and cost of control term with $\nu=1$ in \eqref{eq:tracking-type-Tikhonov}. For the risk measure, an approximation of the conditional value-at-risk measure \eqref{eq:CVAR-def} is used as in \cite{Kouri2019}. In place of the nonsmooth term $v(s) = (1-\beta)^{-1}\max(s,0)$ appearing in the definition of the conditional value-at-risk, the following smooth approximation is used:
\begin{equation*}
\label{eq:smoothed-CVaR}
    v_{{\varepsilon}}(s) = \begin{cases}
    -\frac{\varepsilon}{2}, & \text{if } s\leq -\varepsilon\\
    \frac{1}{2 \varepsilon} s^2 + s, & \text{if } s \in \left(-\varepsilon, \frac{\varepsilon\beta}{1-\beta}\right)\\
    \frac{1}{1-\beta}\left(s-\frac{\varepsilon\beta^2}{2(1-\beta)} \right), & \text{if } s \geq \frac{\varepsilon\beta}{1-\beta}
    \end{cases}
\end{equation*}
with $\varepsilon=0.05$. Note that the smoothed CVaR is still a convex risk measure. The constraint set is given by $\mathbf{K} = \{ v \in V: v \geq 0\}$; i.e., $\psi \equiv 0$. Note that in the numerical section, the state should be greater than or equal to the obstacle. To fit the framework presented in the previous sections, the problem can be transformed using the substitution $\tilde{y}=-y.$

We construct a modification of Example 5.1 from \cite{MR2822818}, an example for which lack of strict complementarity holds (i.e., the measure of the set $\{y^* =0\}\cap\{\zeta^*=0\}$ is positive). We use $A = -\Delta$, {the compact embedding $B\colon L^2(D) \rightarrow H^{-1}(D)$}, and the deterministic functions 
\begin{align*}
    \hat{u}(x)=\hat{y}(x) &= \begin{cases}
    160(x_1^3-x_1^2+0.25x_1)(x_2^3-x_2^2+0.25x_2) & \text{in } (0,0.5)^2,\\
    0 & \text{else},
        \end{cases}\\
    \hat{\zeta}(x)&=  \max(0,-2|x_1-0.8|-2|x_1 x_2-0.3|+0.5),
\end{align*}
constructed in \cite{MR2822818}. Random noise is added to the right-hand side in the form of the (truncated) random field $b\colon D \times \Omega \rightarrow \R$, which is defined by a Karhunen--Lo\`eve expansion; this is described in more detail below. Each random field depends on finite dimensional vectors $\xi\colon \Omega \rightarrow \Xi \subset \R^m$. With these functions, we define the random field $f$ and the target $y_d$ by
\begin{align*}
    f(\cdot,\omega) &:=-\Delta \hat{y}-\hat{y}-\hat{\zeta}-b(\cdot,\omega),\\
    y_d &:= \hat{y}+\hat{\zeta}- \Delta \hat{y}.
\end{align*}
For simulations, problem \eqref{eq:OCProblem} is replaced by the penalised problem \eqref{eq:penalisedOCProblem}; i.e., the inequality constraint is penalised as in \eqref{eq:penalisedPDE} using the smoothed max function defined in \eqref{eq:mrhoHK}. We replace \eqref{eq:penalisedOCProblem} by its sample average approximation (SAA) with the finite set $\Xi = \{\xi_1, \dots, \xi_n\} \subset \R^m$ of randomly drawn vectors. To simplify notation, a sample vector will be denoted by its inverse, i.e., $\omega_i:=\xi_i^{-1}(\omega)$. In summary, the following SAA problems are solved with a decreasing sequence of penalisation parameters $\{\tau_j\}$:
\begin{equation}
\label{eq:SAA-problem}\tag{$P_{\tau}$}
\begin{aligned}
 &\min_{(z,s) \in L^2(D)\times \R} \left \lbrace s+ \frac 1n \sum_{i=1}^n v_\varepsilon\left(\frac 12\norm{y(\cdot,\omega_i)-y_d}{L^2(D)}^2-s \right) + \frac{1}{2}\norm{z}{L^2(D)}^2 \right\rbrace\\
 &\quad \text{s.t.} \quad -\Delta y(\cdot, \omega_\ell) +  \frac 1\tau m_\tau(y(\cdot, \omega_\ell)) =f(\cdot, \omega_\ell)+ z(\cdot), \quad 
\ell=1, \dots, n.
\end{aligned}
\end{equation}
Due to the special structure of CVaR, the control variable is extended by one dimension with $u:=(z,s)$ and $U_{ad}:=L^2(D)\times \R$. {We note that asymptotic consistency (as $n\rightarrow \infty)$ for SAA problems of the form \eqref{eq:SAA-problem} was recently established in \cite[Section 4.3]{Milz2023}, provided the underlying probability space is complete and nonatomic.}

Now we specify our choices for the random field $b$. We observe two examples: one such that $b$ has a pointwise mean zero (in $D$), and the other where $b$ is modelled as a lognormal random field. Both are modifications of examples of random fields on $(-1/2,1/2)^2$ from \cite[{Section 7.4, Example 9.37}]{Lord2014}. These are translated to $D$ and are defined in such a way so that noise is added to only a subset of the biactive set. Examples of realisations of these random fields are displayed in Figure \ref{fig:b}.
\begin{xmpl}[Mean-zero noise]
For the first example, we choose 
\begin{equation*}
\label{eq:truncated-KL-expansion}
 b(x,\omega) = \begin{cases}
\sum_{i=1}^{20} \sqrt{\lambda_i} \phi_i(x) \xi_i(\omega) & \text{in } (0,1/2)\times(0,1),\\
 0 & \text{elsewhere},
 \end{cases}
\end{equation*}
where $\xi_i \sim U(-0.2,0.2)$ for $i=1, \dots, 20$. The eigenfunctions and eigenvalues are given for $j,k\geq 1$ by 
$\tilde{\phi}_{j,k}(x):= 2\cos(j \pi x_2)\cos(k \pi x_1)$ and 
$\tilde{\lambda}_{k,j}:=\frac{1}{4} \exp(-\tfrac{\pi}{4}(j^2+k^2)),$
where we reorder terms so that the eigenvalues appear in descending order (i.e., $\phi_1 = \tilde{\phi}_{1,1}$ and $\lambda_1 = \tilde{\lambda}_{1,1}$).
\end{xmpl}

\begin{xmpl}[{Truncated}  lognormal noise]
In this example, noise is added to a subset of the biactive set in the form of lognormal field with truncated Gaussian noise by
\begin{equation*}
b(x,\omega) = \begin{cases}
e^{-4+\sum_{i=1}^{100} \sqrt{\lambda_i} \phi_i(x) \xi_i(\omega)} & \text{in } (0,1/2)\times(0,1/2),\\
0 & \text{elsewhere},
\end{cases}
\end{equation*}
where $\xi_i$ distributed according to the truncated normal distribution $\mathcal{N}(0,3,-100,100)$ with mean $0$ and standard deviation $3.$ The eigenfunctions $\phi_j(x) = \phi_{i,1}(x_1)\phi_{k,2}(x_2)$ and eigenvalues $\lambda_j = \lambda_{i,1} \lambda_{k,2}$ are given by the following functions (relabeled after sorting by decreasing eigenvalues):
\begin{equation*}
\begin{aligned}
 \phi^{i,m}(x_m) &= \begin{cases}
               \sqrt{1/2+\sin(w_i)/(2 w_i)}^{-1} \cos(w_i x_m) & \text{for } i \text{ odd},\\ 
                \sqrt{1/2-\sin(w_i)/(2 w_i)}^{-1} \sin(w_i x_m) & \text{for } i \text{ even},
              \end{cases}\\
\lambda_{i,m} &= \frac{2}{w_i^2 + 1}, \quad \quad \quad  w_i = \begin{cases}
             \hat{w}_{\lceil{i/2}\rceil} & \text{for } i \text{ odd},\\
             \tilde{w}_{i/2} & \text{for } i \text{ even},
            \end{cases}
\end{aligned}
\end{equation*}
where $\hat{w}_j$ is the $j^{\text{th}}$ positive root of $1-w \tan(w/2)$, and $\tilde{w}_j$ is the $j^{\text{th}}$ positive root of $\tan(w/2) + w.$ 
\begin{figure}
    \centering
    \includegraphics[height=4cm]{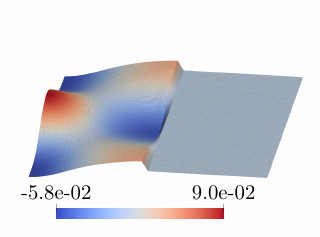}
 \includegraphics[height=4cm]{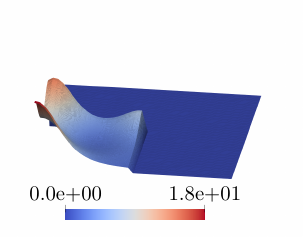}
\includegraphics[height=4cm]{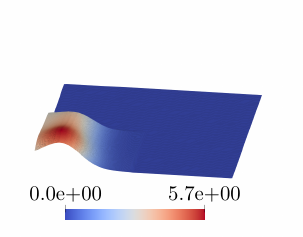}
    \caption{Example realisations of random field zero-mean $b$ (left) and lognormal random fields $b$ (middle, right).}
    \label{fig:b}
\end{figure}
\end{xmpl}

\subsection{Path-following stochastic approximation}
In numerical experiments, we solve a sequence of the SAA-approximated proxy problems \eqref{eq:SAA-problem} and iteratively decrease the penalisation term $\tau$ in an outer loop. The proxy problems are solved using the stochastic variance reduced gradient (SVRG) method from \cite{johnson2013accelerating}. {The method reduces variance by regularly computing a full gradient for the SAA problem \eqref{eq:SAA-problem}, balanced by intermediate random steps. This computational strategy has several advantages for our problem: first, intermediate gradients do not need to be stored (unlike other methods like stochastic average gradient). Second, the regular computation of the full gradient provides a certificate of optimality in order to terminate the inner loop (SVRG) procedure. Finally, the computational cost is greatly reduced when compared to standard batch gradient methods; see \cite{johnson2013accelerating}.} An alternative to using this low iteration complexity method would be to apply the primal dual risk minimization method from \cite{kouri2022primal} to the problem \eqref{eq:SAA-problem}. The advantage of this approach would be higher accuracy for each $\tau$ at a possibly higher computational cost.
Let
\begin{equation*}
    J(u,\omega):=s + v_\varepsilon\left(\frac{1}{2}\norm{y(\cdot,\omega)-y_d}{L^2(D)}^2-s \right) + \frac{1}{2}\norm{z}{L^2(D)}^2
\end{equation*}
be the parametrised objective function corresponding to the problem \eqref{eq:SAA-problem}. For the algorithm, we rely on a stochastic gradient $G_{{\tau}}\colon L^2(D)\times \Omega \rightarrow L^2(D)$, i.e., the function satisfying $\E[G_{{\tau}}(u,\cdot)] = \nabla \E[J(u,\cdot)].$ 
The stochastic gradient is defined by 
\begin{equation*}
 G_{{\tau}}(u,\omega) = 
\begin{pmatrix}
 p_{{\tau}}(\cdot,\omega)+ z\\
 1-v_\varepsilon'(\frac{1}{2} \norm{y_{{\tau}}(\cdot, \omega)-y_d}{L^2(D)}^2 - s)
  \end{pmatrix},
\end{equation*}
where $p_{{\tau}}(\cdot, \omega)$ solves \eqref{eq:ptau} and $y_{{\tau}}(\cdot, \omega)$ solves  \eqref{eq:ytau}.
The full gradient for the SAA approximation is denoted by 
$g_{{\tau}}(u)=\frac{1}{n}\sum_{i=1}^n G_{{\tau}}(u,\omega_i).$
For the termination of the middle loop, we use the residual 
\begin{equation*}
\hat{r}_{{\tau}}(u):= \norm{\frac{1}{n}\sum_{i=1}^n p_{{\tau}}(\cdot,\omega_i)+ z}{L^2(D)}+\left\vert 1-v_\varepsilon'\left(\frac{1}{2} \norm{\frac{1}{n}\sum_{i=1}^n y_{{\tau}}(\cdot,\omega_i)-y_d}{L^2(D)}^2 - s\right)\right\vert.
\end{equation*}

\begin{algorithm}[H] 
\begin{algorithmic}[1] 
\STATE \textbf{Initialisation:} Choose $\tilde{u}_1$, penalty smoothing parameter $\tau_1$, update frequency $r$, step-size sequence $\{ t_{k_\ell} \}$, tolerance \texttt{tol}, smoothing multiplier $\gamma$, $k=1$
\FOR{$j=1,2,\dots$}
\WHILE{$\hat{r}_{{\tau_j}}(u) > \texttt{tol}$} 
\STATE $u_1 :=\tilde{u}_k$
\STATE $\hat{g} := g_{{\tau_j}}(u_1)$
\STATE Randomly sample $n_k$ from $\{ 1, \dots, r\}$ 
\FOR{$\ell=1, 2, \dots, n_k$}
\STATE Randomly sample $i_\ell$ from $\{1, \dots, n\}$
\STATE Set $u_{\ell+1} := u_\ell - t_{k_\ell}( G_{{\tau_j}}(u_{\ell},\omega_{i_\ell})-G_{{\tau_j}}(u_1,\omega_{i_\ell})+\hat{g})$
\ENDFOR
\STATE $k:=k+1$
\STATE $\tilde{u}_k:=u_{n_k+1}$
\ENDWHILE
\STATE $\tau_{j+1}:=\gamma\tau_j$
\ENDFOR
\end{algorithmic}
\captionof{algorithm}{Path-following SVRG} 
\label{alg:SVRG}
\end{algorithm}
In the simulations, the random indices from lines 7 and 9 in Algorithm \ref{alg:SVRG} are generated according to the uniform distribution, i.e., $n_k \sim \mathcal{U}(\{1, \dots, r \})$ and $i_\ell \sim \mathcal{U}(\{1, \dots, n\})$, although other choices are possible. While the convergence of Algorithm \ref{alg:SVRG} has yet to be proven in the function space setting, the convergence of the stochastic gradient method (without variance reduction) has been established; see \cite{Geiersbach2021d, Geiersbach2020b} for convergence results when the method is applied to nonconvex PDE-constrained optimisation problems.

All simulations were done using Python along with the finite element environment FEniCS (2018.1.0) \cite{Alnes2015}. For the generation of random numbers, we use \texttt{numpy.random.seed(4)}. This seed is used to generate random numbers in the following order: first, for the $i_k$, a random vector of length 5,000 is generated (according to the discrete uniform distribution over $\{1, \dots, r\}$ with update frequency $r=1,000$). Then, 10,000 vectors of length $m$ ($m=20$ for the zero-mean example and $m=100$ for the lognormal example) are produced. Finally, for each $k$, a random vector of length $i_k$ is generated (according to the discrete uniform distribution over $\{1, \dots, n\}$.) Full gradient computations are done with the help of the \texttt{multiprocessing} module. All functions are discretised using P2 Lagrange finite elements with {3200 elements}.  The state equation \eqref{eq:ytau} is solved with relative tolerance $10^{-8}$ using a Newton solver. The adjoint equation
\eqref{eq:ptau} is solved with a relative tolerance $10^{-8}$ using the Krylov solver GMRES with the ILU preconditioner.

Termination conditions are informed by \cite{MR2822818}. The tolerance in the middle loop is chosen to be \texttt{tol} = {$6.25\cdot 10^{-7}$}. 
The smoothing multiplier in the outer loop is chosen to be $\gamma=0.1$. The step-size from the original SVRG method \cite{johnson2013accelerating} is a constant that depends on the Lipschitz constant of the gradient and strong convexity, which are clearly not available. 
Consequently, we use the step-size rule 
\begin{equation*}
t_{k_\ell}= \sqrt{\frac{\theta}{k \ell+\nu}}, \quad \theta = \frac{1}{2\nu}+1, \quad \nu = \frac{2 \theta}{2\nu-1}-1
\end{equation*}
inspired by a similar rule developed for convex problems in \cite{Geiersbach2020b}. 

For starting values, we choose $u_1\equiv 1$, and $s_1 = 1$. We remark that a proper choice of $s_1$ appears to greatly impact the performance of the method. This value was chosen to be in the neighborhood of the first $s_k$
such that second component of the gradient satisfies $|1-v_\varepsilon'(\frac{1}{2} \norm{\E[y]-y_d}{L^2(D)}^2 - s)| <1.$ 

\begin{figure}[hbtp]
    \centering
        \includegraphics[height=3.4cm]{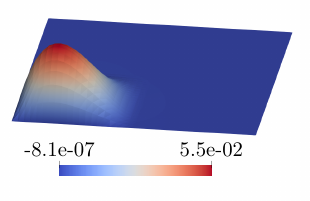}
 \includegraphics[height=3.4cm]{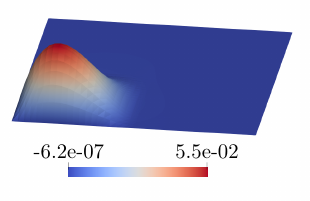}
\includegraphics[height=3.4cm]{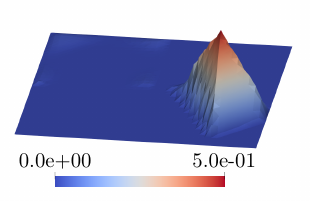}
    \caption{Mean-zero noise: Control $u_\tau^*$ (left), state $\bar{y}_\tau^*$ (middle), multiplier $\bar{\zeta}_\tau^*$ (right) for $\beta=0.0$.}
    \label{fig:solutions_zero_0}

\vspace*{1\floatsep}
\includegraphics[height=4cm]{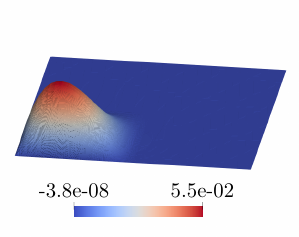}
 \includegraphics[height=4cm]{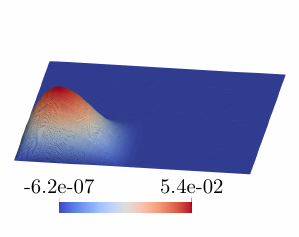}
\includegraphics[height=4cm]{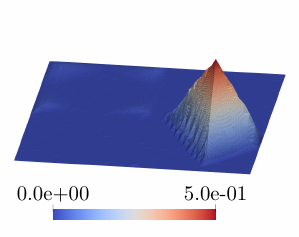}
    \caption{Mean-zero noise: Control $u_\tau^*$ (left), state $\bar{y}_\tau^*$ (middle), multiplier $\bar{\zeta}_\tau^*$ (right) for $\beta=0.95$.}
    \label{fig:solutions_zero_95}
\vspace*{3\floatsep}
        \includegraphics[height=3.4cm]{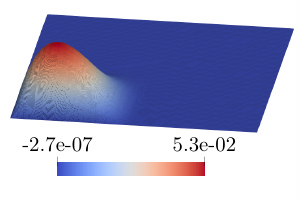}
 \includegraphics[height=3.4cm]{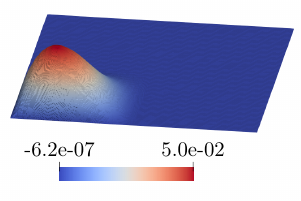}
\includegraphics[height=3.4cm]{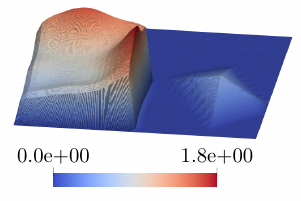}
    \caption{Lognormal noise: Control $u_\tau^*$ (left), state $\bar{y}_\tau^*$ (middle), multiplier $\bar{\zeta}_\tau^*$ (right) for $\beta=0.0$.}
    \label{fig:solutions_log_0}
\vspace*{1\floatsep}
        \includegraphics[height=4cm]{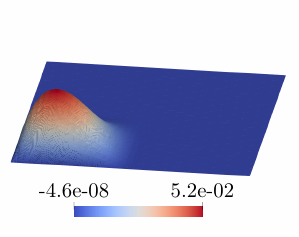}
 \includegraphics[height=4cm]{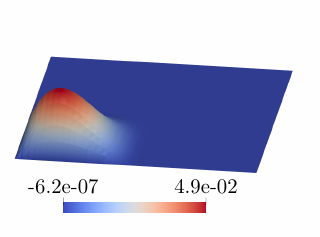}
\includegraphics[height=4cm]{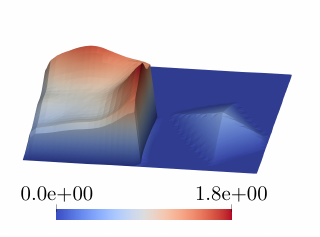}
    \caption{Lognormal noise: Control $u_\tau^*$ (left), state $\bar{y}_\tau^*$ (middle), multiplier $\bar{\zeta}_\tau^*$ (right) for $\beta=0.95$.}
    \label{fig:solutions_log_95}
\end{figure}

In Table \ref{tab:numerical_values}, numerical values are displayed showing the final objective value obtained for $\tau=10^{-6}$ with the aforementioned error tolerance. {A reviewer of this paper noted that the risk-neutral and CVaR solutions being very similar is indicative of a problem with low variance.} The risk-averse examples are significantly more expensive than their risk-neutral counterparts and the lognormal case, which exhibits higher variance than the mean-zero case and requires more PDE solves for the given risk level $\beta$.
\begin{table}[!htb]
\begin{center}
 \begin{tabular}{ c|  c | c | c | c }
&  Mean-zero, $\beta=0$ & Mean-zero, $\beta=0.95$ &  Lognormal, $\beta=0$ & Lognormal, $\beta=0.95$ \\ 
  \hline
$\bar{j}_\tau^*$ &  1.3952056 & 1.3954905 & 1.3965497 & 1.3969087\\
\# Full gradient comp. & 14 & 37 & 16 & 48\\
\# PDE solves & 291,808 & 758,108 & 342,104 & 1,063,756
\end{tabular}
 \end{center}
 \caption{Final objective function value $\bar{j}_\tau^*$ achieved for $\tau=10^{-6}$ and computational cost}
 \label{tab:numerical_values}
\end{table}
The control $u_\tau$ and averaged solutions
\begin{equation*}
\bar{y}_\tau^* = \frac{1}{n} \sum_{i=1}^n y_{{\tau}}(\cdot, \omega_i), \quad \bar{\zeta}_\tau^* = \frac{1}{n} \sum_{i=1}^n \zeta_{{\tau}}(\cdot, \omega_i)
\end{equation*}
are shown for different levels of CVaR$_{\beta}$ in Figures \ref{fig:solutions_zero_0} to \ref{fig:solutions_log_95} for $\tau=10^{-6}$. One sees here that the lack of strict complementarity persists in the averaged solutions.
In the case of mean-zero noise, the solutions resemble the deterministic {solutions (shown in Figure~\ref{fig:solutions_deterministic})} as expected.
The risk-averse case with $\beta=0.95$ shows a slight difference in the minimal and maximal values of the solutions and states. For the lognormal case, where the variance of the random field is also greater, we see greater differences in function values in Figures \ref{fig:solutions_log_0} to \ref{fig:solutions_log_95}. The additional noise is above all apparent in the multiplier $\bar{\zeta}_\tau.$

\begin{figure}
    \centering
        \includegraphics[height=3.6cm]{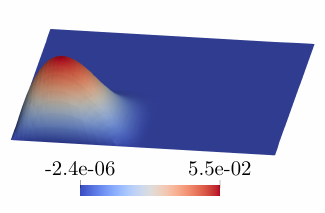}  
         \label{fig:solutions_deterministic}
        \includegraphics[height=3.6cm]{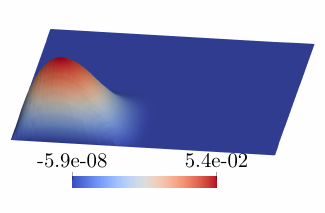}  
        \includegraphics[height=3.6cm]{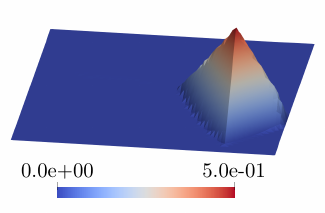}  
    \caption{{Deterministic solution: Control $u_\tau^*$ (left), state $\bar{y}_\tau^*$ (middle), multiplier $\bar{\zeta}_\tau^*$ (right).}}
      \label{fig:solutions_deterministic}
\end{figure}

\section{Conclusion}\label{sec:conclusion}

In this paper, we focused on VIs of obstacle type in an $L^2(D)$ setting. We could instead have worked in an abstract Gelfand triple setting $(V,H,V^*)$ and with a more general assumption on the constraint set for the VI \eqref{eq:VIProblemFormal} such as
\[\text{$\mathbf{K}(\omega) \subset V$ is a non-empty, closed and convex subset}.\]
The maps $m_\tau$ would need to be modified for this abstract setting, see \cite[\S 2.3]{AHROCQVI}. Results up to and including the existence of optimal controls should hold with the same assumptions on the various operators (with obvious modifications where necessary) but the derivation of stationarity conditions would require more thought (see the comments after Theorem 5.5 of \cite{AHROCQVI}). 

{Let us also reiterate that the stationarity systems that we obtained in this paper are indeed weaker than their counterparts in the deterministic setting, see the discussions after Theorem \ref{thm:ss-weaker} and Proposition \ref{thm:ss}. It is an interesting open problem to resolve these missing items.}

One could also attempt to tackle the nonsmooth problem \eqref{eq:OCProblem} directly (making use of the VI directional differentiability results of \cite{MR0423155}) without the penalisation approach we took here, but then it is rather unclear how to unfold the primal conditions and obtain a dual stationarity system.

While the proposed path-following SVRG method was able to compute solutions up to a high accuracy, the method is not yet optimised and theoretical justification on the appropriate function spaces is missing. In particular, development of step-size rules coupled with a mesh refinement strategy in the style of \cite{Geiersbach2020b} will be the topic of future research. 
As already mentioned, the primal dual algorithm for risk minimization \cite{kouri2022primal} could have been used to solve the smoothed subproblems; a future study could include a rigorous hypothesis test as in \cite{kouri2023relaxation} to compare the application of a first-order method like SVRG to the primal dual algorithm. {Once we opt for smoothed CVaR in the extended setting $U_{ad} \times \mathbb R$, the optimization problem fits into the general risk-neutral setting considered in \cite{hertlein2022inexact}. Thus, in principle, it is possible to proceed numerically with the infinite-dimensional bundle method analyzed there. However, the method in \cite{hertlein2022inexact} is for the fully continuous setting, which suggests that an empirical approximation/SAA should be employed to treat the objective and remain in their setting. From a practical side, this would result in an objective function defined as the sum of potentially thousands of nonsmooth, nonconvex functions over an infinite dimensional space. An interesting future direction would be to investigate the possibility of using stochastic subgradients in the context of the bundle method in \cite{hertlein2022inexact} in order to reduce the computational burden. As a side note, the only known mesh independent methods for solving the obstacle problem use some form of penalty or barrier method, e.g., \cite{Hintermller2006}, which would mean that the evaluation of $S$ in the nonsmooth setting is either to be done with a mesh dependent method, e.g., the primal dual active set strategy and semismooth Newton methods \cite{Hintermller2002,Ulbrich2002}, or that separate penalty inner loops for each sample VI are required.}

\appendix
\section{Differentiability of superposition operators}
Combining Theorem 7 and Remark 4 of \cite{Goldberg}\footnote{Note that $p$ and $q$ are switched in that paper to what we have here.}, we have the following lemma.

 \begin{lmm}\label{lem:goldberg}
Let $\alpha$ be a number satisfying $1 \leq p < \alpha < \infty$  and let $X$ and $Y$ be (real) Banach spaces. Suppose $H \colon X \times \Omega \to Y$ is a Carath\'eodory function that is Fr\'echet differentiable with respect to $x \in X$ and assume that $H_x\colon X \times \Omega \to \mathcal{L}(X,Y)$ is a Carath\'eodory function. Furthermore, 
     assume there exists $C_1 \in L^p(\Omega)$ and $C_2 \geq 0$ such that
     \begin{equation*}
        \norm{H(x,\omega)}{Y} \leq C_1(\omega) + C_2\norm{x}{X}^{\alpha\slash p} \quad \text{a.s.} \quad \forall x \in X
    \end{equation*}
and assume that there exists $\tilde{C}_1 \in L^{\tilde r}(\Omega)$ where $\tilde r = p\alpha\slash(\alpha-p)$ and $\tilde{C}_2 \geq 0$ such that
\begin{equation*}
\norm{H_x(x,\omega)}{\mathcal{L}(X,Y)} \leq \tilde{C}_1(\omega) + \tilde{C}_2\norm{x}{X}^{\alpha/\tilde r} \quad \text{a.s.} \quad \forall x \in X.
\end{equation*}
{Then the Nemytskii operator $\mathcal{H}$ defined through the formula
\[(\mathcal{H}u)(\omega) = H(u(\omega),\omega)\quad \text{for $u\colon \Omega \to X$}\]
is such that} $\mathcal{H}\colon L^\alpha(\Omega;X) \to L^p(\Omega; Y)$ is continuously Fr\'echet differentiable {and the derivative $\mathcal{H}' \colon L^\alpha(\Omega;X) \rightarrow \mathcal{L}(L^\alpha(\Omega;X), L^p(\Omega;Y))$ is given by}
\[(\mathcal{H}'(u)h)(\omega)  = H_x(u(\omega), \omega)h(\omega) \quad\text{for $\omega\in \Omega$, $u, h \in L^\alpha(\Omega;X)$}.\]
In addition, { 
$\mathcal{H}'\colon L^\alpha(\Omega;X) \to L^{p\alpha\slash (\alpha-p)}(\Omega;\mathcal{L}(X,Y))$ is continuous.}
 \end{lmm}
 \begin{proof}
 The conditions on $H$ imply that $\mathcal{H}$ maps $L^\alpha(\Omega;X)$ to $L^p(\Omega;Y)$ due to \cite[Theorem 1]{Goldberg}. Under the conditions on $H_x$, we find that the Nemytskii operator  $\hat{\mathcal{H}}' \colon L^\alpha(\Omega;X) \to L^{s}(\Omega; \mathcal{L}(X,Y))$  is a continuous map \cite[Theorems 4 and 5]{Goldberg} (see also \cite[Remark 4]{Goldberg}) where $s = p\alpha\slash (\alpha-p)$ in the first case.  
 Then we simply apply \cite[Theorem 7]{Goldberg}.
 \end{proof}
 Note that the Fr\'echet differentiability in combination with $\alpha < p \leq \infty$ would imply $\mathcal{H}$ is constant, whereas with $p=\alpha < \infty$ implies affineness and thus these cases are excluded from the above.
\section{Other results}\label{sec:pqStuff}
{Let  $p,q \in (1,\infty)$ be exponents and define
\[\tilde p:= \frac{pq}{q-p}.\] 
Recall that we use the notation $p'$ to denote the conjugate of $p$.}
\begin{lmm}\label{lem:tildePGeqTwo}
Let $p<q.$ If 
\begin{equation*}
\text{either $p \geq 2$ or if $2 \leq q \leq 2p$,}    
\end{equation*}
we have $\tilde p \geq 2$.
\end{lmm}
\begin{proof}
If $p \geq 2$, this follows immediately from \eqref{eq:tildePGreaterThanP} and in fact in this case we get strict inequality. Consider
\begin{align*}
    \frac{pq}{q-p}-2 &= \frac{pq}{q-p} - \frac{2q-2p}{q-p}= \frac{pq-2q+2p}{q-p}\geq \frac{4p-2q}{q-p}
\end{align*}
and this is non-negative if $q \leq 2p.$

\end{proof}
\begin{lmm}
Let $p<q$. If $q \leq 2p$, then $\tilde p \geq q$.
\end{lmm}
\begin{proof}

Consider
\begin{align*}
    \frac{pq}{q-p}-q &= \frac{pq}{q-p} - \frac{q^2-pq}{q-p}= \frac{2pq-q^2}{q-p}= \frac{q(2p-q)}{q-p}
\end{align*}
and this is non-negative if $q \leq 2p.$
\end{proof}

\begin{lmm}\label{lem:trivialLemma}
If $\tilde p \geq 2$, we have $\tilde p \geq \tilde p'.$
\end{lmm}
\begin{proof}
We have
\begin{align*}
    \tilde p - \tilde p' &= \tilde p - \frac{\tilde p}{\tilde p-1}= \frac{\tilde p^2 - 2\tilde p}{\tilde p -1}= \frac{\tilde p(\tilde p - 2)}{\tilde p -1}
\end{align*}
which is non-negative when $\tilde p \geq 2$.
\end{proof}
\bibliographystyle{abbrv}
\bibliography{references}
\end{document}